\documentclass[11pt]{article}
\pdfoutput =1 
\usepackage[left=1in,right=1in,top=1in,bottom=1in]{geometry}
\oddsidemargin-2mm\evensidemargin-2mm\usepackage{amsmath}
\usepackage{color}
\usepackage[normalem]{ulem}
\usepackage{xcolor}
\newcommand\redsout{\bgroup\markoverwith{\textcolor{red}{\rule[0.5ex]{2pt}{0.4pt}}}\ULon}
\newcommand{\Limsup}{\mathop{{\rm Lim}\,{\rm sup}}}
\def\disp{\displaystyle}

\def\tto{\;{\lower 1pt\hbox{$\rightarrow$}}\kern-10pt
\hbox{\raise 2pt\hbox{$\rightarrow$}}\;}
\def\Hat{\widehat}
\def\hat{\widehat}

\def\tilde{\widetilde}
\def\Bar{\overline}
\def\ra{\rangle}
\def\la{\langle}
\def\ve{\varepsilon}
\def\B{I\!\!B}
\def\IN{I\!\!N}

\def\h{\hfill\Box}
\def\R{I\!\!R}

\def\K{\mathcal{K}}
\def\ox{\bar{x}}
\def\op{\bar{p}}
\def\oy{\bar{y}}

\def\ov{\bar{v}}

\def\oq{\bar{q}}
\def\ow{\bar{w}}

\def\gph{\mbox{\rm gph}\,}

\def\epi{\mbox{\rm epi}\,}
\def\dim{\mbox{\rm dim}\,}
\def\dom{\mbox{\rm dom}\,}

\def\h{\hfill\triangle}
\def\dn{\downarrow}
\def\O{\Omega}
\def\Lm{\Lambda}
\def\ph{\varphi}
\def\emp{\emptyset}
\def\st{\stackrel}
\def\oR{\Bar{\R}}

\def\lm{\lambda}
\def\gg{\gamma}
\def\dd{\delta}
\def\al{\alpha}
\def\kk{\kappa}

\def\Th{\Theta}
\def\th{\theta}
\def\vt{\vartheta}
\newcounter{lk}
\def\Limsup{\mathop{{\rm Lim}\,{\rm sup}}}

\begin{document}
\newtheorem{Theorem}{Theorem}[section]
\newtheorem{Proposition}[Theorem]{Proposition}
\newtheorem{Remark}[Theorem]{Remark}
\newtheorem{Lemma}[Theorem]{Lemma}
\newtheorem{Corollary}[Theorem]{Corollary}
\newtheorem{Definition}[Theorem]{Definition}
\newtheorem{Example}[Theorem]{Example}
\renewcommand{\theequation}{\thesection.\arabic{equation}}
\normalsize
\def\proof{
\normalfont
\medskip
{\noindent\itshape Proof.\hspace*{6pt}\ignorespaces}}
\def\endproof{$\h$\vspace*{0.05in}}
\begin{center}
\vspace*{0.3in} {\bf FULL STABILITY OF GENERAL PARAMETRIC VARIATIONAL SYSTEMS}\\[2ex]

B. S. MORDUKHOVICH\footnote{Department of Mathematics, Wayne State University, Detroit, MI 48202, USA and the RUDN University, Moscow 117198, Russia (boris@math.wayne.edu). Research of this author was partly supported by the US National Science Foundation under grant DMS-1512846, by the US Air Force Office of Scientific Research under grant \#15RT0462, and by the Ministry of Education and Science of the Russian Federation (Agreement number 02.a03.21.0008 of 24 June 2016).}, T. T. A. NGHIA\footnote{Department of Mathematics and Statistics, Oakland University, Rochester, MI 48309, USA (nttran@oakland.edu).} and D. T. PHAM \footnote{Department of Mathematics, Wayne State University, Detroit, MI 48202, USA (dat.pham@wayne.edu). Research of this author was partly supported by the US National Science Foundation under grant DMS-1512846 and by the US Air Force Office of Scientific Research under grant \#15RT0462.}
\end{center}

\small{\bf Abstract.} The paper introduces and studies the notions of Lipschitzian and H\"olderian full stability of solutions to three-parametric variational systems described in the generalized equation formalism involving nonsmooth base mappings and partial subgradients of prox-regular functions acting in Hilbert spaces. Employing advanced tools and techniques of second-order variational analysis allows us to establish complete characterizations of, as well as directly verifiable sufficient conditions for, such full stability notions under mild assumptions. Furthermore, we derive exact formulas and effective quantitative estimates for the corresponding moduli. The obtained results are specified for important classes of variational inequalities and variational conditions in both finite and infinite dimensions.\\[1ex]
\vspace*{-0.15in}

\normalsize	
\section{Introduction}
\vspace*{-0.05in}

This paper concerns the area of {\em second-order variational analysis}, which has been of a rapidly increasing interest during the recent years. Our main attention is paid to the study of {\em full stability} for general {\em three-parametric variational systems} (PVS) given by
\begin{equation}\label{1.1}
v\in f(x,p,q)+\partial_x g(x,p),
\end{equation}
where $x\in X$ stands for the {\em decision} variable from a Hilbert space $X$, where $(v,p,q)\in X\times\mathcal{P}\times\mathcal{Q}$ is a triple of perturbation {\em parameters} with $v\in X$ signifying {\em canonical} perturbations while $(p,q)\in\mathcal{P}\times\mathcal{Q}$ is the pair of {\em basic} perturbations taking values in metric spaces $\mathcal{P}$ and $\mathcal{Q}$, where $f\colon X\times\mathcal{P}\times\mathcal{Q}\to X$ is a single-valued {\em base} mapping, and where $g\colon X\times\mathcal{P}\to\oR:=(-\infty,\infty]$ is an extended-real-valued and lower semicontinuous (l.s.c.) {\em potential} function for which the symbol $\partial_x$ indicates the set of its {\em partial limiting subgradients} with respect to the decision variable; see Section~2 for more details.

Some particular cases of the PVS model \eqref{1.1}, including variational and quasi-variational inequalities, variational conditions, etc., have been intensively studied by many researchers over the years with numerous applications to nonlinear analysis, optimization, equilibria, ordinary and partial differential equations, control theory, and numerical algorithms; see, e.g., \cite{DR,FP,hms,KS,M1,R0,R03,rw} and the bibliographies therein as well as the further references presented below. The pioneering albeit fundamental impacts on these directions of research and applications were done by Stampacchia \cite{S} for infinite-dimensional variational inequalities motivated by applications to PDEs, and by Robinson \cite{R0} in the framework of generalized equations motivated mainly by applications to numerical optimization. The major thrust in both publications was on obtaining efficient conditions ensuring the properties of {\em single-valuedness} and {\em continuity} or {\em Lipschitz continuity} of the corresponding solution maps associated with such important specifications of PVS.

In this paper we thoughtfully investigate the notions of Lipschitzian and H\"olderian full stability for PVS \eqref{1.1} involving behavior of their parameter-dependent {\em solution maps}
\begin{equation*}
S(v,p,q):=\big\{x\in X\big|\;v\in f(x,p,q)+\partial_x g(x,p)\big\},\quad(v,p,q)\in X\times\mathcal{P}\times\mathcal{Q},
\end{equation*}
with respect to perturbations of the reference parameter triple $(\ov,\op,\oq)$. Such notions have been recently introduced by Mordukhovich and Nghia \cite{MN2} and studied therein in the case where the base $f$ does not depend on the parameter $q$ while being smooth with respect to the state variable $x$. These full stability notions provide nontrivial extensions to PVS of the corresponding definitions for local minimizers originated by Levy, Poliquin and Rockafellar \cite{LPR} in the Lipschitzian case and then proceeded in \cite{MN} for the case of H\"olderian full stability.

Considering here the more general version \eqref{1.1} allows us to unify various frameworks of PVS studied in the literature being important for a variety of applications. It is not hard to show that the results of \cite{MN2} can be readily extended to the general PVS framework of \eqref{1.1} provided that the base mapping $f$ is smooth in $x$. Thus the goals of this paper are different from such an extension.

Our {\em first goal} is to derive efficient conditions ensuring the H\"olderian and Lipschitzian full stability for nonsmooth PVS \eqref{1.1} with replacing the smoothness of $f$ by a certain local {\em strong monotonicity} of this mapping with respect to the state variable. Furthermore, the full stability of such PVS is justified under a {\em quantitative} relationship between the modulus of strong monotonicity of $f$ and the threshold of prox-regularity of the potential $g$ in \eqref{1.1}, which is a new result even for smooth PVS with the base $f$ independent of $q$ as in \cite{MN2}.

The {\em second goal} is motivated by the above while being definitely important for its own sake. We obtain, for the first time in the literature, {\em exact formulas} for calculating the {\em threshold of prox-regularity} of $g$ in \eqref{1.1} in terms of the {\em second-order subdifferential} (generalized Hessian) constructions initiated by the first author \cite{M92}. Besides the usage in this paper, the obtained calculations can be utilized in the development and justification of numerical algorithms of optimization; in particular, of the proximal point and related types; see, e.g., \cite{P} and the references therein.

The {\em third goal} of the paper is to provide implementations and specifications of the general results obtained on full stability of \eqref{1.1} for particular classes of PVS that are well-recognized in variational theory and applications. There are two classes of such systems considered in the paper: {\em variational inequalities} in Hilbert spaces and the so-called {\em variational conditions} in finite dimensions. For the first class we derive {\em pointbased characterizations} of Lipschitzian full stability under some {\em polyhedrality} conditions. The second class is described by the normal cone mappings generated by {\em parameter-dependent} sets given by smooth inequalities. The obtained second-order characterizations of full stability are expressed entirely via the initial data of such systems.\vspace*{0.03in}

We organize the rest of the paper as follows. Section~2 recalls some {\em preliminaries} of variational analysis and generalized differentiation that are widely employed below. We also formulate here a generalized notion of uniform second-order growth condition in the vein of Bonnans and Shapiro \cite{BS}. The next Section~3 presents the {\em basic definitions} of H\"olderian and Lipschitzian full stability of the three-parametric PVS \eqref{1.1} and then lists and discusses our {\em standing assumptions} for the most of the paper. In Section~4 give extensions of some results from \cite{MN2} to the case of \eqref{1.1} under the partial {\em smoothness} of $f$ with respect to $x$.

Section~5 is devoted to the qualitative and quantitative study of parametrically continuous prox-regular functions $g\colon X\times\mathcal{P}\to\oR$ that appear as potentials in \eqref{1.1}. We provide such a study via second-order subdifferential constructions for $g$ and establish, in particular, exact formulas for computing the threshold of prox-regularity in the second-order subdifferential terms.

Sections 6 and 7 concern deriving verifiable conditions, involving the interplay between the strong monotonicity modulus of the nonsmooth base $f$ and the threshold of prox-regularity for the potential $g$ in \eqref{1.1}, that ensure the H\"olderian and Lipschitzian full stability of PVS. Our approach here is completely different from \cite{MN2} while relaying on the study of the corresponding properties for the nonconvex extension of the so-called {\em proximal mapping} for $g$ well-recognized in optimization theory and algorithms. The obtained results for the proximal mapping certainly are of their independent interest. We present discussions and examples on the relationships the new results with related ones in the literature; in particular, with those established by Yen \cite{Y2,Y1}.

Section~8 deals with {\em parametric variational inequalities} (PVI) written as
\begin{eqnarray}\label{vi}
v\in f(x,p)+N_C(x)\;\mbox{ for }\;x\in C\subset X,\;p\in \mathcal{P},
\end{eqnarray}
where $N_C(x)$ is the normal cone at $x$ to the closed {\em convex} subset $C$ of the Hilbert space $X$. It is clear that \eqref{vi} is a particular case of PVS \eqref{1.1} with $g(x)=\delta_C(x)$ being the indicator function of the parameter-independent convex set $C$. It follows from the normal cone construction for convex set that \eqref{vi} can be rewritten as follows: given $p\in\mathcal{P}$, find $x\in C$ such that
\begin{eqnarray}\label{vi1}
\la v-f(x,p),u-x\ra\le 0\;\mbox{ for all }\;u\in C.
\end{eqnarray}
Infinite-dimensional variational inequalities in form \eqref{vi1} often appear in optimization-related (in particular, optimal control) problems governed by elliptic partial differential equations, which are usually modeled via the so-called {\em Legendre form} under the {\em polyhedricity} assumption on $C$; see, e.g., \cite{BBS,B,BS,hms,HS,IT,KS} and the precise definitions in Section~7 for more details. Imposing these natural requirements and elaborating the results in Section~3, we derive in Section~7 {\em pointbased} (i.e., expressed exactly at the reference points) {\em characterizations} of Lipschitzian full stability of solutions to the perturbed variational inequalities \eqref{vi} held in the infinite-dimensional framework of Hilbert decision spaces in spite of the {\em lack of compactness} in infinite dimensions.

In Section~9 we study another type of PVS \eqref{1.1} that are known as {\em parametric variational conditions} (PVC) and are represented by
\begin{eqnarray}\label{qvi}
v\in f(x,p,q)+N_{C(p)}(x)\;\mbox{ with }\;x\in C(p)\subset X,\;p\in\mathcal{P},\;q\in\mathcal{Q}
\end{eqnarray}
via the limiting normal cone $N_{C(p)}$ to $C(p)$ at $x$, where the decision space $X$ and the parameter spaces $\mathcal{P,Q}$ are finite-dimensional, and where the parameter-dependent set $C(p)$ is described by the finitely many inequality constraints
\begin{eqnarray}\label{ine}
C(p):=\big\{x\in X\big|\;\ph_i(x,p)\le 0\;\mbox{ for }\;i=1,\ldots,m\big\}
\end{eqnarray}
defined by ${\cal C}^2$-smooth functions. Such systems are clearly different from \eqref{vi} due to the dependence of the moving set $C(\cdot)$ on the perturbation parameter $p$; they reduce to parametric {\em quasi-variational inequalities} if the sets $C(p)$ are convex; see, e.g., \cite{FP,K2,L2,LR,mo,R03,R13,Y2,Y1} and the references therein for various terminology and stability results concerning parametric systems of type \eqref{qvi}.

We introduce in Section~9 a new second-order qualification condition under the name of ``general uniform second-order sufficient condition" (GUSOSC) and show that it completely characterizes Lipschitzian full stability of solutions to \eqref{qvi} under the simultaneous validity of the partial Mangasarian-Fromovitz and constant rank constraint qualifications (MFCQ and CRCQ, respectively) for \eqref{ine}. If both these constraint qualifications are replaced by the stronger linear independence constraint qualification (LICQ) for the partial gradients of the active constraints in \eqref{ine} at the reference point, then the new GUSOSC reduces to the known ``general strong second-order sufficient condition" (GSSOSC) from \cite{K2}, a slight modification for variational conditions of Robinson's strong second-order sufficient condition \cite{R} in parametric nonlinear programming with ${\cal C}^2$-smooth data. In this way we arrive at a rather surprising result that GSSOSC is sufficient for Lipschitzian full stability of \eqref{qvi} and {\em completely characterizes} this property under LICQ. We present an example showing that our GUSOSC holds and thus ensures, in particular, the properties of local single-valuedness and Lipschitz continuity for the solution map to \eqref{qvi} with linear constraints in \eqref{ine} while the well-recognized ``strong coherent orientation condition" known to be sufficient for these properties in such a setting \cite{FP} fails to fulfill.\vspace*{0.03in}

Throughout of the paper we use the standard notation and terminology of variational analysis; cf. \cite{M1,rw}. Unless otherwise stated, the {\em decision} space $X$ is {\em Hilbert} while being identified with its dual space $X^*$. Recall that $\la\cdot,\cdot\ra$ stands for the canonical pairing in $X$ with the norm $\|x\|:=\sqrt{\la x,x\ra}$ and $\st{w}{\to}$ denotes the weak convergence in $X$. Furthermore, $\B$ indicates the closed unit ball in the space in question, and thus $\B_\eta(x):=x+\eta\B$ is the closed ball centered at $x$ with radius $\eta>0$. For a set-valued mapping $F\colon X\tto X$ from $X$ into itself $(X=X^*)$, the symbol
\begin{equation}\label{pk}
\begin{array}{ll}
\disp\Limsup_{x\to\ox}F(x):=\Big\{v\in X\Big|&\exists\;\mbox{ sequences }\;x_k\to\ox,\;v_k\st{w}{\to}v\;\mbox{ such that}\\
&v_k\in F(x_k)\;\mbox{ for all }\;k\in\IN:=\{1,2,\ldots\}\Big\}
\end{array}
\end{equation}
signifies the (sequential) {\em Painlev\'e-Kuratowski outer/upper limit} of $F(x)$ as $x\to\ox$. As stated at the beginning, the {\em parameter} space $(\mathcal{P},d_1)$ is {\em metric} with $\B_\eta(p):=\{p'\in P|\;d_1(p',p)\le\eta\}$ standing for the closed ball centered at $p$ with radius $\eta>0$. Similar notations are used in the metric space $(\mathcal{Q},d_2)$. The closed ball in the product space $X\times\mathcal{P}$ is referred as $\B_\eta(x,p):=\B_\eta(x)\times\B_\eta(p)$. Finally, that the symbols $x\st{\ph}{\to}\ox$ and $x\st{\O}{\to}\ox$ for a function $\ph\colon X\to\oR$ and a set $\O\subset X$, respectively, mean that $x\to\ox$ with $\ph(x)\to\ph(\ox)$ and $x\in\O$.\vspace*{-0.15in}

\section{Preliminaries from Variational Analysis}	
\setcounter{equation}{0}\vspace*{-0.05in}

We begin with reviewing the generalized differential constructions broadly used in the paper and refer the reader to \cite{M1,rw} for more details. Given an extended-real-valued function $\ph\colon X\to\oR$ with $\ox$ from $\dom\ph:=\{x\in X|\;\ph(x)<\infty\}$, the {\em regular subdifferential} of $\ph$ at $\ox$
(known also as the presubdifferential as well as the Fr\'echet or viscosity subdifferential of $\ph$ at $\ox$) is given by
\begin{eqnarray}\label{2.1}
\Hat\partial\ph(\ox):=\Big\{v\in X\Big|\;\liminf_{x\to\ox}\frac{\ph(x)-\ph(\ox)-\la v,x-\ox\ra}{\|x-\ox\|}\ge 0\Big\}.
\end{eqnarray}
The {\em limiting subdifferential} of $\ph$ at $\ox$ (known also as the basic or Mordukhovich subdifferential) and its {\em singular subdifferential} counterpart are defined, respectively, via the outer limit \eqref{pk} by
\begin{eqnarray}\label{2.2}
\partial\ph(\ox):=\Limsup_{x\st{\ph}{\to}\ox}\Hat\partial\ph(x)\;\mbox{ and }\;\partial^\infty\ph(\ox):=\Limsup_{x\st{\ph}{\to}\ox,\lm\dn 0}\lm\Hat\partial\ph(x).
\end{eqnarray}
When $\ph$ is convex, both regular and limiting subdifferentials above reduce to the subdifferential of convex analysis while $\partial^\infty\ph(\ox)=\{0\}$ if $\ph$ is locally Lipschitzian around $\ox$. It is worth mentioning that, despite (perhaps due to) the nonconvexity of the subgradient sets in \eqref{2.2}, these subdifferentials and the corresponding normal cone and coderivative constructions for sets and mappings posses {\em full calculi} derived from the {\em variational/extremal principles} of variational analysis; see \cite{M1,rw}.

The corresponding {\em regular} and {\em limiting normal cones} to $\O$ at $\ox\in\O$ are defined by
\begin{eqnarray}\label{nc}
\Hat N_\O(\ox):=\Hat\partial\delta_\O(\ox)\;\mbox{ and }\;N_\O(\ox):=\partial\delta_\O(\ox)
\end{eqnarray}
via the subdifferential constructions \eqref{2.1} and \eqref{2.2} applied to the indicator function $\delta_\O(x)$ of $\O$ equal to $0$ for $x\in\O$ and to $\infty$ otherwise.

Consider further a set-valued mapping $F\colon X\tto Y$ between two Hilbert spaces with its domain and graph defined in the standard way by
$$
\dom F:=\big\{x\in X\big|\;F(x)\ne\emp\big\}\;\mbox{ and }\;\gph F:=\big\{(x,y)\in X\times Y\big|\;y\in F(x)\big\}.
$$
The {\em regular} and {\em limiting coderivatives} of $F$ at $(\ox,\oy)$ are given, respectively, as
\begin{eqnarray}
\Hat D^*F(\ox,\oy)(w)&:=&\big\{z\in X\big|\;(z,-w)\in\Hat N_{{\rm gph} F}(\ox,\oy)\big\},\quad w\in Y,\label{2.3}\\
D^*F(\ox,\oy)(w)&:=&\big\{z\in X\big|\;(z,-w)\in N_{{\rm gph}F}(\ox,\oy)\big\},\quad w\in Y,\label{2.4}
\end{eqnarray}
where we skip $\oy=F(\ox)$ from the coderivative notation when $F$ is single-valued around $\ox$.

A mapping $F\colon X\tto Y$ from a metric space $(X,d)$ to a normed space $Y$ is {\em Lipschitz-like} (pseudo-Lipschitz or having the Aubin property) around $(\ox,\oy)\in\gph F$ with modulus $\ell>0$ if
\begin{equation}\label{2.5}
F(x)\cap V\subset F(u)+\ell d(x,u)\B\;\mbox{ for all }\;x,u\in U
\end{equation}
for some neighborhoods $U$ of $\ox$ and $V$ of $\oy$. As discussed in \cite[Chapter~1]{M1}, the inclusion \eqref{2.5} is equivalent to the distance estimate
\begin{equation}\label{2.6}
{\rm haus}\big(F(x)\cap V,F(u)\big)\le\ell d(x,u)\;\mbox{ for all }\;x,u\in U,
\end{equation}
where haus$(C_1,C_2)$ is the Pompieu-Hausdorff distance between the sets $C_1$ and $C_2$ defined by
$$
{\rm haus}(C_1,C_2):=\inf\big\{\eta\ge 0\big|\;C_1\subset C_2+\eta\B,\,\;C_2\subset C_1+\eta\B\big\}.
$$
When both $X$ and $Y$ are finite-dimensional spaces, the Lipschitz-like property is fully characterized by the following injectivity condition:
\begin{eqnarray}\label{cod-cr}
D^*F(\ox,\oy)(0)=\{0\}
\end{eqnarray}
known as the coderivative/Mordukhovich criterion; see \cite[Corollary~5.4]{m93} and \cite[Theorem~9.40]{rw}. The infinite-dimensional extensions of this result can be found in \cite[Theorems~4.7 and 4.10]{M1}.\vspace*{0.03in}

Next we present the standard versions of single-valued localizations of set-valued mappings as defined, e.g., in \cite{DR,MN2} and used throughout the paper.\vspace*{-0.1in}

\begin{Definition}{\bf (single-valued localizations).}\label{locali} Let $F\colon X\tto Y$ be a set-valued mapping between metric spaces, and let $(\ox,\oy)\in \gph F$. We say that $F$ admits a {\sc single-valued localization} around $(\ox,\oy)$ if there is a neighborhood $U\times V\subset X\times Y$ of $(\ox,\oy)$ such that the mapping $\Hat F\colon U\to V$ defined by $\gph\Hat F:=\gph F\cap(U\times V)$ is single-valued on $U$ with $\dom\Hat F=U$. In this case we say that $\Hat F$ is a single-valued localization of $F$ relative to $U\times V$. If in addition $\Hat F$ is $($Lipschitz$)$ continuous on $U$, then $F$ admits a {\sc $($Lipschitz$)$ continuous single-valued localization} around $(\ox,\oy)$.
\end{Definition}\vspace*{-0.05in}

Finally in this section, we recall the uniform second-order growth condition for extended-real-valued functions, which was formulated in \cite[Definition~3.5]{mrs} while reducing to the so-called ``uniform second-order (quadratic) growth condition with respect to the ${\cal C}^2$-smooth parameterization" for problems of ${\cal C}^2$ conic programming introduced in \cite[Definition~5.16]{BS}.\vspace*{-0.1in}

\begin{Definition}{\bf (uniform second-order growth condition).}\label{USOGC} Let $X$ and $\mathcal{P}$ be a Hilbert and metric space, respectively. Given $h\colon X\times\mathcal{P}\to\oR$ and $\ov\in\partial_x h(\ox,\op)$, we say the {\sc uniform second-order growth condition (USOGC)} holds at $(\ox,\op,\ov)$ with modulus $\ell>0$ if there are neighborhoods $U$ of $\ox$, $V$ of $\ov$, and $P$ of $\op$ such that
\begin{eqnarray}\label{2.10}
h(x,p)\ge h(u,p)+\la v,x-u\ra+\frac{\ell}{2}\|x-u\|^2\;\mbox{ for }\;x\in U,\;(u,p,v)\in\gph\partial_{x}h\cap(U\times P\times V).
\end{eqnarray}
\end{Definition}\vspace*{-0.3in}

\section{Basic Notions and Standing Assumptions}
\setcounter{equation}{0}\vspace*{-0.05in}

Consider a general parametric variational system given in the form
\begin{eqnarray}\label{VS}
v\in f(x,p,q)+\partial_x g(x,p)\;\mbox{ for }\;x\in X,\;p\in\mathcal{P},\;q\in\mathcal{Q}
\end{eqnarray}
with the Hilbert decision space $X$ and the metric parameter spaces $(\mathcal{P},d_1)$ and $(\mathcal{Q},d_2)$, where $f\colon X\times \mathcal{P}\times\mathcal{Q}\to X$, $g\colon X\times\mathcal{P}\to\oR$, and $\partial_x g$ stands for the partial limiting subdifferential of the function $g$ with respect to the variable $x$ taken from \eqref{2.2}. Denote $g_p(\cdot):=g(\cdot,p)$ and observe that $\partial_x g(x,p)=\partial g_p(x)$ for all $(x,p)\in X\times\mathcal{P}$. Fixing $\ov\in f(\ox,\op,\oq)+\partial_x g(\ox,\op)$ and define further the solution map $S\colon X\times\mathcal{P}\times\mathcal{Q}\tto X$ to PVS \eqref{VS} by
\begin{eqnarray}\label{ss}
S(v,p,q):=\big\{x\in X\big|\;v\in f(x,p,q)+\partial_x g(x,p)\big\}\;\mbox{ with }\;\ox\in S(\ov,\op,\oq).
\end{eqnarray}
Such a general formalism of variational analysis was investigated in \cite{MN2} in the case of $f=f(x,p)$ therein with the base mapping $f$ smooth in $x$, while the definitions of H\"olderian and Lipschitzian full stability for PVS introduced therein do not depend on the aforementioned smoothness. The next formulations of full stability for \eqref{VS} adapt the scheme of \cite{MN2} with the {\em different treatment} of the parameters $p$ and $q$ in the H\"olderian version.\vspace*{-0.1in}

\begin{Definition}{\bf (H\"olderian and Lipschitzian full stability for general PVS).}\label{fs} Given a solution $\ox\in S(\ov,\op,\oq)$ to PVS \eqref{VS} corresponding to the parameter triple $(\ov,\op,\oq)$, we say that:
	
{\bf (i)} $\ox$ is a {\sc H\"olderian fully stable} solution to \eqref{VS} corresponding to $(\ov,\op,\oq)$ if the solution map \eqref{ss} has a single-valued localization $\vt$ relative to some neighborhood $V\times P\times Q\times U$ of $(\ov,\op,\oq,\ox)$ so that for any $(v_1,p_1,q_1), (v_2,p_2,q_2)\in V\times P\times Q$ there are constants $\kk,\ell>0$ with
\begin{eqnarray}\label{4.7}
\big\|(v_1-v_2)-2\kk[\vt(v_1,p_1,q_1)-\vt(v_2,p_2,q_2)]\big\|\le\|v_1-v_2\|+\ell\big[d_1(p_1,p_2)^{\frac{1}{2}}+d_2(q_1,q_2)\big].
\end{eqnarray}
	
{\bf (ii)} $\ox$ is a {\sc Lipschitzian fully stable} solution to \eqref{VS} corresponding to $(\ov,\op,\oq)$ if the solution map \eqref{ss} has a single-valued localization $\vt$ relative to some neighborhood $V\times P\times Q\times U$ of $(\ov,\op,\oq,\ox)$ so that for any $(v_1,p_1,q_1),(v_2,p_2,q_2)\in V\times P\times Q$ there are constants $\kk,\ell>0$ with
\begin{eqnarray}\label{mp1}
\big\|(v_1-v_2)-2\kk\big[\vt(v_1,p_1,q_1)-\vt(v_2,p_2,q_2)\big]\big\|\le\|v_1-v_2\|+\ell\big[d_1(p_1,p_2)+d_2(q_1,q_2)\big].
\end{eqnarray}
\end{Definition}\vspace*{-0.05in}

From the first glance the full stability notions for PVS from Definition~\ref{VS} have nothing to do with the corresponding notions of Lipschitzian \cite{LPR} and H\"olderian \cite{MN} full stability of local minimizers. However, it follows from \cite[Corollary~4.6]{MN2} that in the case of $f=0$ in \eqref{VS} the full stability notions in Definition~\ref{VS} agree with those for local minimizers of the potential function $g$ in \eqref{VS}. This provides {\em new error bound characterizations} \eqref{4.7} and \eqref{mp1} of full stability for local minimizers. Furthermore, it is worth emphasizing that the stability notions from Definition~\ref{fs} clearly yield the {\em local single-valuedness} and {\em H\"older} (resp. {\em Lipschitz$)$ continuity} of the solution map \eqref{ss}, but not vice versa as the simple examples in \cite{MN} demonstrate.\vspace*{0.03in}

Before formulating the standing assumptions on the data of PVS \eqref{VS} imposed in this paper, we recall the definition of a major parametric class of extended-real-valued functions that play a fundamental role in second-order variational analysis and its applications. The parametric version used in what follows was defined in \cite{LPR}, while the original nonparametric notion was given in \cite{PR1}; see also \cite{BT1,ST,rw} and the references therein for more details in finite and infinite dimensions. Let $g\colon X\times\mathcal{P}\to\oR$ be l.s.c.\ around $(\ox,\op)\in\dom g$ with $\hat v:=\ov-f(\ox,\op,\oq)\in\partial_x g(\ox,\op)$. Following \cite{LPR}, we say that $g$ is {\em prox-regular} in $x$ at $\ox$ for $\hat v$ with {\em compatible parameterization} by $p$ at $\op$ if there exist neighborhoods $U$ of $\ox$, $V$ of $\hat v$, and $P$ of $\op$ along with positive numbers $\ve,r$ satisfying
\begin{eqnarray}\label{par-prox}
\begin{array}{ll}
\quad\quad g(x,p)\ge g(u,p)+\la v,x-u\ra-\frac{r}{2}\|x-u\|^2\;\mbox{ for all }\;x\in U,\\
\quad\mbox{when }\;v\in\partial_x g(u,p)\cap V,\;u\in U,\;p\in P,\;\mbox{ and }\;g(u,p)\le g(\ox,\op)+\ve.
\end{array}
\end{eqnarray}
The infimum of all $r$ in \eqref{par-prox} is called the {\em threshold of prox-regularity} of $g$ at $(\ox,\op)$ for $\hat v$ and is denoted by $\mathcal{R}$. The function $g$ is said to be {\em subdifferentially continuous} in $x$ at $\ox$ for $\hat v$ with {\em compatible parameterization} by $p$ at $\op$ if the mapping $(x,p,v)\mapsto g(x,p)$ is continuous relative to the set $\gph\partial_x g$ at $(\ox,\op,\hat v)$. For brevity $g$ is called {\em parametrically continuously prox-regular} at $(\ox,\op)$ for $\hat v$ if it is prox-regular and subdifferentially continuous at $\ox$ for $\hat v$ with compatible parameterization by $p$ at $\op$. In the latter case the inequality ``$g(u,p)\le g(\ox,\op)+\ve$" is extra in \eqref{par-prox}.\vspace*{0.03in}

Now we are ready to formulate the {\em standing assumptions} used throughout the paper:\vspace*{0.03in}

{\bf (A1)} The base $f$ is {\em Lipschitz continuous} around $(\ox,\op,\oq)$, i.e., there are a constant $L>0$ and a neighborhood $U\times P\times Q$ of $(\ox,\op,\oq)$ such that for all $ x_1,x_2\in U,\,p_1,p_2\in P,\;q_1,q_2\in Q$ we have
\begin{equation}\label{lipf}
\|f(x_1,p_1,q_1)-f(x_2,p_2,q_2)\|\le L\big[\|x_1-x_2\|+d_1(p_1,p_2)+d_2(q_1,q_2)\big].
\end{equation}

{\bf (A2)} The potential $g$ is {\em parametrically continuously prox-regular} at $(\ox,\op)$ for $\hat v$.

{\bf (A3)} The {\em basic constraint qualification} (BCQ) is satisfied at $(\ox,\op)$:
\begin{eqnarray}\label{bcq1}
\mbox{the mapping}\;\;p\mapsto\epi g(\cdot,p)\;\;\mbox{is {\em Lipschitz-like} around }\;\big(\op,(\ox,g(\ox,\op))\big).
\end{eqnarray}

If $\dim X\cdot\dim {\cal P}<\infty$, assumption (A3) can be equivalently reformulated via the singular subdifferential in \eqref{2.2} as follows
\begin{eqnarray}\label{bcq2}
(0,p)\in\partial^\infty g(\ox,\op)\Longrightarrow p=0,
\end{eqnarray}
which is an immediate consequence of the coderivative criterion \eqref{cod-cr}. As shown in Sections~8 and 9, assumptions (A2) and (A3) are automatically satisfied for important special classes of parametric variational systems in both finite and infinite dimensions. Note, in particular, that for the indicator function  $g(x,p)=\delta_{C(p)}(x)$ of a closed-valued multifunction $C\colon P\tto X$ the BCQ assumption (A3) is equivalent to the Lipschitz-like property of $C$ around $(\op,\ox)$ used in \cite{Y2,Y1} for parametric variational inequalities. Indeed, we have in this case that $\epi g(\cdot,p)=C(p)\times\R^+$. Hence
\begin{eqnarray}\label{haus}
\begin{array}{ll}
{\rm haus}\big(\epi g(\cdot,p_1)\cap(V\times\R),\epi g(\cdot,p_2)\big)&\disp={\rm haus}\big((C(p_1)\cap V)\times\R^+,C(p_2)\times\R^+\big)\\
&\disp={\rm haus}\big(C(p_1)\cap V,C(p_2)\big)
\end{array}
\end{eqnarray}
for a neighborhood $V$ of the reference point, and the aforementioned equivalence follows from \eqref{2.6}.\vspace*{-0.15in}

\section{Full Stability of PVS under Partial Differentiability}
\setcounter{equation}{0}\vspace*{-0.05in}

To derive second-order characterizations and verifiable sufficient conditions for H\"olderian and Lipschitzian full stability of PVS, we use throughout the paper
the second-order subdifferential constructions for extended-real-valued functions defined by the scheme \cite{M92} as a {\em coderivative} of their {\em first-order subdifferentials}. Dealing here with the parametric potential function $g_p(x):=g(x,p)$ in \eqref{VS} and its partial subdifferential $\partial_x g(x,p)=\partial g_p(x)$ requires for our purposes the usage the following two modifications of the {\em second-order partial subdifferentials} defined by
\begin{eqnarray}\label{sec}
(D^*\partial g_p)(x,p,v)=(D^*\partial_x g)(x,p,v)\;\mbox{ and }\;(\Hat D^*\partial g_p)(x,p,v)=(\Hat D^*\partial_x g)(x,p,v)
\end{eqnarray}
defined via the limiting and regular coderivatives, respectively, of the partial subdifferential
\begin{eqnarray*}
\partial_x g(x,p)=\big\{\mbox{{\rm set of limiting subgradients of }}\;g_p=g(\cdot,p)\;\mbox{ {\rm at }}\;x\big\}.
\end{eqnarray*}
The reader can find more information on various properties and applications of these partial second-order constructions in \cite{LPR}, \cite{mnn}, \cite{MN}, \cite{MN2}, \cite{MNR13}, \cite{MR}, \cite{mrs}, and the references there. The usage of the regular vs.\ limiting coderivative in \eqref{sec} was first suggested in \cite{MN13} in the nonparametric setting.

This section addresses second-order characterizations of H\"olderian and Lipschitzian full stability of \eqref{VS} under the additional partial smoothness assumption on the base mapping:\vspace*{0.03in}

{\bf (A4)} $f$ is (Fr\'echet) differentiable with respect to $x$ around $(\ox,\op,\oq)$ uniformly in $(p,q)$, and the partial Jacobian $\nabla_x f$ is continuous at $(\ox,\op,\oq)$.\vspace*{0.03in}

The results given below under (A4) are extensions of those obtained in \cite{MN2} in the case of $f=f(x,p)$ in \eqref{VS}, i.e., when both base and potential in \eqref{VS} depend on the same parameter.

First we present characterizations of Lipschitzian full stability in \eqref{VS} while starting with a {\em neighborhood} conditions at points near the reference one obtained in infinite dimensions.\vspace*{-0.1in}

\begin{Proposition}{\bf(neighborhood second-order characterization of Lipschitzian full stability of general PVS).}\label{Lippq} Let $\ox\in S(\ov,\op,\oq)$ be a solution to \eqref{VS} corresponding to $(\ov,\op,\oq)$, and let assumptions {\rm(A1)--(A4)} be satisfied. Consider the following conditions:
	
{\bf (i)} There exist positive numbers $\eta,\kk_0$ such that whenever $(u,p,v)\in\gph\partial_x g\cap\B_\eta(\ox,\op,\hat v)$ with $\hat v=\ov-f(\ox,\op,\oq)$ we have the condition
\begin{eqnarray}\label{4.8b}
\la\nabla_x f(\ox,\op,\oq)w,w\ra+\la z,w\ra\ge\kk_0\|w\|^2\quad\mbox{for all}\quad z\in(\Hat D^*\partial g_p)(u,v)(w),\;w\in X.
\end{eqnarray}
	
{\bf (ii)} The graphical set-valued partial subdifferential mapping
\begin{eqnarray}\label{K}
K\colon p\mapsto\gph\partial_x g(\cdot,p)\;\;\mbox{is Lipschitz-like around}\;\;(\op,\ox,\hat v).
\end{eqnarray}
Then the validity of both conditions {\rm(i)} and {\rm(ii)} is equivalent to the Lipschitzian full stability of $\ox$.
\end{Proposition}\vspace*{-0.07in}
{\bf Proof.} Suppose that both conditions (i) and (ii) are satisfied. Define the auxiliary extended-real-valued function $\tilde g\colon X\times\mathcal{P}\times\mathcal{Q}\to\oR$ by $\tilde g(x,p,q):=g(x,p)$ and observe that condition (ii) of the proposition is equivalent to the Lipschitz-like property of the auxiliary graphical mapping $\tilde{K}\colon (p,q)\mapsto\gph\partial_x\tilde g(\cdot,p,q)$ around $(\op,\oq,\ox,\hat v)$. Define now the extended PVS
\begin{eqnarray}\label{4.10}
v\in f\big(x,(p,q)\big)+\partial_x\tilde g\big(x,(p,q)\big),
\end{eqnarray}
which is equivalent to \eqref{VS}. We can treat \eqref{4.10} as a PVS depending on the common parameter $\tilde p=(p,q)$ in both single-valued and set-valued parts of the variational system. Applying now \cite[Theorem~4.7]{MN2} to the one-parameter PVS \eqref{4.10} and taking into account the constructions of the auxiliary data of the latter system verify the claimed Lipschitzian full stability of \eqref{VS}. The converse implication of the proposition can be checked similarly by reducing it to \cite[Theorem~4.7]{MN2}.\endproof

Note that in the case of $\dim X,\dim{\cal P}<\infty$ the Lipschitz-like requirement in \eqref{K} can be completely characterized via the coderivative criterion \eqref{cod-cr} by the pointbased second-order condition
\begin{eqnarray}\label{Mor}
(0,z)\in\big(D^*\partial_x g\big)(\ox,\op,\hat v)(0)\Longrightarrow z=0.
\end{eqnarray}

The next result gives a {\em pointbased} characterization of Lipschitzian full stability for \eqref{VS} in finite dimensions via the limiting partial second-order subdifferential from \eqref{sec}.\vspace*{-0.1in}

\begin{Proposition}{\bf (pointbased characterization of Lipschitzian full stability of general PVS).}\label{Lipspq} Let $X,P$ be finite-dimensional in the framework of Proposition~{\rm\ref{Lippq}}. Then  $\ox$ is a Lipschitzian fully stable in \eqref{VS} if and only if condition \eqref{Mor} holds simultaneously with
\begin{eqnarray}\label{4.43}
\la\nabla_x f(\ox,\op,\oq)w,w\ra+\la z,w\ra>0\;\mbox{ for all }\;(z,s)\in\big(D^*\partial_x g\big)(\ox,\op,\hat v)(w),\;w\ne 0.
\end{eqnarray}
\end{Proposition}\vspace*{-0.05in}
{\bf Proof.} It goes in the same way as the proof of Proposition~\ref{Lippq} with applying now the result of \cite[Theorem~4.8]{MN2} for one-parameter PVS instead of \cite[Theorem~4.7]{MN2} in the proof above.\endproof

The situation with H\"olderian full stability for PVS \eqref{VS} is more involved in comparison with the Lipschitzian case. In this setting we arrive at the following result.\vspace*{-0.1in}

\begin{Theorem}{\bf (second-order characterization of H\"olderian full stability for general PVS).}\label{Holderpq} Let $\ox\in S(\ov,\op,\oq)$ in \eqref{ss} under assumptions {\rm(A1)--(A4)}. Consider the assertions:
	
{\bf (i)} The solution $\ox$ is a H\"olderian fully stable in \eqref{VS} with the moduli $\kk,\ell>0$ from \eqref{4.7} corresponding to the parameter triple $(\ov,\op,\oq)$.
	
{\bf (ii)} For some positive numbers $\eta,\kk_0$ and any $(u,p,v)\in\gph\partial_x g\cap\B_\eta(\ox,\op,\hat v)$ with $\hat v=\ov- f(\ox,\op,\oq)$ we have the condition
\begin{eqnarray}\label{4.8a}
\la\nabla_x f(\ox,\op,\oq)w,w\ra+\la z,w\ra\ge\kk_0\|w\|^2\quad\mbox{whenever}\quad z\in(\Hat D^*\partial g_p)(u,v)(w),\;w\in X.
\end{eqnarray}
Then it holds {\rm(i)}$\Longrightarrow${\rm(ii)} with constant $\kk_0$ that can be chosen smaller than but arbitrarily close to $\kk$. Conversely we get {\rm(ii)}$\Longrightarrow${\rm(i)}, where $\kk$ can be chosen smaller but arbitrarily close to $\kk_0$.
\end{Theorem}\vspace*{-0.07in}
{\bf Proof.} Assuming that (i) is satisfied and using Definition~\ref{fs}(i), we find all the parameters of H\"olderian full stability satisfying inequality \eqref{4.7}. Pick $\eta\in(0,1)$ sufficiently small. Suppose without loss of generality that $Q=\B_{\frac{\eta}{2}}(\oq)$ and that
$$
d_2(q_1,q_2)\le\sqrt{\eta d_2(q_1,q_2)}\;\mbox{ for any }\;q_1,q_2\in Q.
$$
Then it clearly follows from \eqref{4.7} that the inequality
\begin{equation}\label{4.7a}
\big\|(v_1-v_2)-2\kk[\vt(v_1,p_1,q_1)-\vt(v_2,p_2,q_2)]\big\|\le\|v_1-v_2\|+\ell\big[d_1(p_1,p_2)^{\frac{1}{2}}+ d_2(q_1,q_2)^{\frac{1}{2}}\big]
\end{equation}
holds for any $(v_1,p_1,q_1),(v_2,p_2,q_2)\in V\times P\times Q$. Note that
$$
2d_{\mathcal{P}\times\mathcal{Q}}((p_1,q_1),(p_2,q_2))=2 \big[d_1(p_1,p_2)+d_2(q_1,q_2)\big]\ge\big[d_1(p_1,p_2)^{\frac{1}{2}}+d_2(q_1,q_2)^{\frac{1}{2}}\big]^2.
$$
This together with \eqref{4.7a}  gives us the condition
\begin{equation}\label{4.7c}
\big\|(v_1-v_2)-2\kk[\vt(v_1,p_1,q_1)-\vt(v_2,p_2,q_2)]\big\|\le\|v_1-v_2\|+\sqrt{2}\ell\, d_{\mathcal{P}\times\mathcal{Q}}((p_1,q_1),(p_2,q_2))^{\frac{1}{2}}.
\end{equation}
Similarly to the proof of Proposition~\ref{Lippq}, we form the extended parameter $\tilde p=(p,q)$ and observe that condition \eqref{4.7c} ensures the H\"olderian full stability of the extended one-parameter PVS \eqref{4.10}. Employing now implication (i)$\Longrightarrow$(ii) in \cite[Theorem~4.3]{MN2} justifies the validity of the corresponding second-order condition for \eqref{4.10}, which is equivalent to \eqref{4.8a}.

To verify next the opposite implication of this theorem, suppose that (ii) holds and conclude from implication (ii)$\Longrightarrow$(i) of \cite[Theorem~4.3]{MN2} that there are numbers $\kk,\ell>0$ and a neighborhood $V\times P\times Q\times U$ of $(\ov,\op,\oq,\ox)$ on which the solution map to the extended system \eqref{4.10} admits a single-valued localization $\vartheta$ satisfying
\begin{equation}\label{4.11}
\big\|(v_1-v_2)-2\kk[\vt(v_1,p_1,q_1)-\vt(v_2,p_2,q_2)]\big\|\le\|v_1-v_2\|+\ell\big[d_{\mathcal{P}\times\mathcal{Q}}((p_1,q_1),(p_2,q_2))\big]^{\frac{1}{2}}
\end{equation}
for all triples $(v_1,p_1,q_1),(v_2,p_2,q_2)\in V\times P\times Q$. Pick any two such triples and let $x_1:=\vt(v_1,p_1,q_1)$, $x_2:=\vt(v_2,p_2,q_2)$, and $x_3:=\vt(v_2,p_1,q_2)$. Then it follows from \eqref{4.11} that
\begin{eqnarray}\label{4.12}
2\kk\|x_2-x_3\|\le\ell\sqrt{d_1(p_1,p_2)}.
\end{eqnarray}
Fixing finally $p_1$ and applying Proposition~\ref{Lippq} to the PVS
$$
v\in f(x,p_1,q)+\partial_x g(x,p_1)
$$
tells us that $\ox$ is Lipschitzian fully stable for this PVS, and hence we get
\begin{eqnarray}\label{4.13}
\|(v_1-v_2)-2\kk[x_1-x_3]\big\|\le\|v_1-v_2\|+\tilde\ell d_2(q_1,q_2)
\end{eqnarray}
with some positive constant $\tilde\ell$. Combining the latter with \eqref{4.12} and using the triangle inequality implies that $\ox$ is a H\"olderian fully stable solution of the PVS \eqref{VS} under consideration.
\endproof
\vspace*{-0.2in}

\section{Calculating the Threshold of Prox-Regularity}
\setcounter{equation}{0}

This section is devoted to the study of the parametrically continuous {\em prox-regularity} of the potential function $g$ in PVS \eqref{VS} that is defined in Section~3. Besides the applications of the results obtained here to the study of full stability of nonsmooth PVS in the subsequent sections, these results are of their own importance for other applications. In particular, we derive exact formulas for computing the {\em threshold} of prox-regularity via our second-order generalized differential constructions.

Let us start with two lemmas one of which is taken from \cite{MN} while the other one constitutes a new result of its independent interest. The first lemma follows directly from the combination of \cite[Theorem~4.5 and Theorem~4.7]{MN}. \vspace*{-0.1in}

\begin{Lemma}{\bf (second-order subdifferential characterization of USOGC).}\label{fullmin} Let $X$ and ${\cal P}$ be a Hilbert and metric space, respectively, and let $h\colon X\times\mathcal{P}\to\R$ be parametrically continuously prox-regular at $(\ox,\op)$ for $\hat{v}\in\partial_x h(\ox,\op)$. Then under the validity of BCQ from assumption {\rm (A3)} for $h$ the following statements are equivalent:

{\bf (i)} USOGC \eqref{2.10} from Definition~{\rm\ref{USOGC}} holds at $(\ox,\op,\hat{v})$ with some modulus $\ell>0$.

{\bf (ii)} There is positive number $\eta$ such that whenever $(u,p,v)\in\gph\partial_x h\cap\B_\eta(\ox,\op,\hat{v})$ we have
$$
\la z,w\ra\ge\ell\|u\|^2\;\mbox{ for all }\;z\in(\Hat D^*\partial h_p)(u,v)(w),\;w\in X,
$$
where $h_p(x):=h(x,p)$ with $\partial_x h(x,p)=\partial h_p(x)$
\end{Lemma}\vspace*{-0.05in}

The next lemma presents a verifiable necessary condition for the prox-regularity under consideration expressed in the second-order subdifferential terms. \vspace*{-0.1in}

\begin{Lemma}{\bf (second-order necessary condition for prox-regularity).}\label{lm0} Let $\ox\in S(\ov,\op,\oq)$ in \eqref{ss} under the validity of the standing assumptions {\rm(A2)} and {\rm(A3)} on the potential $g$, and let ${\cal R}$ be the threshold of prox-regularity in {\rm(A2)}.
Given any $r>\mathcal{R}$, there is a neighborhood $U\times P\times V$ of $(\ox,\op,\hat v)$ with $\hat v=\ov-f(\ox,\op,\oq)$ such that for every $(u,p,v)\in\gph\partial_x g\cap(U\times P\times V)$ we have
\begin{equation}\label{4.14}
\la z,w\ra\ge-r\|w\|^2\quad\mbox{whenever}\quad z\in(\Hat D^*\partial g_p)(u,v)(w),\;w\in X.
\end{equation}
\end{Lemma}
{\bf Proof.} Since $g$ is parametrically continuously prox-regular with the prox-parameter $r$ at $(\ox,\op)$ for $\hat v$, there exist neighborhoods $U$ of $\ox$, $V$ of $\hat v$, and $P$ of $\op$ with
\begin{equation}\label{e.9.1}
g(x,p)-g(u,p)\ge\la v,x-u\ra-\frac{r}{2}\|x-u\|^2\;\;\mbox{for all}\;\;(u,p,v)\in\gph\partial_x g\cap\big(U\times P\times V\big),\;x\in U.
\end{equation}
Picking any $s>r$, define $h(x,p):=g(x,p)+\frac{s}{2}\|x-\ox\|^2$ and use the notation $h_p$ as above. We are going to show that $h$ satisfies the conditions USOGC and BCQ from \eqref{2.10} and \eqref{bcq1}, respectively. To proceed with USOGC first, deduce from the subdifferential sum in \cite[Proposition~1.107]{M1} that $\partial h_p(\ox)=\partial g_p(\ox)+s(I-\ox)$ for any $p\in\mathcal{P}$. Define now $W:=J(U\times P\times V )$ with $J(x,p,v):=(x,p,v+s(x-\ox))$ and conclude from the open mapping theorem that $W$ is a neighborhood of $(\ox,\op,\hat v)$. It is easy to observe the validity of the inclusion
$$
v-s(u-\ox)\in\partial g_p(u)\cap V\;\mbox{ for any }\;(u,p,v)\in\gph\partial_x h\cap W.
$$
This together with \eqref{e.9.1} tells us that for any $x\in U$ we have
\begin{eqnarray*}
\begin{array}{ll}
h_p(x)&\disp=g_p(x)+\frac{s}{2}\|x-\ox\|^2\\
&\disp\ge g_p(u)+\la v-s(u-\ox),x-u\ra-\frac{r}{2}\|x-u\|^2+\frac{s}{2}\|x-\ox\|^2\\
&\disp=h_p(u)+\la v,x-u\ra+\frac{s-r}{2}\|x-u\|^2.
\end{array}
\end{eqnarray*}
As a result, it shows that USOGC \eqref{2.10} holds for $h$ at $(\ox,\op)$ with modulus $(s-r)$.

We can easily see that $h$ is parametrically continuously prox-regular at $(\ox,\op)$ for $\hat v$ and proceed now with verifying that BCQ \eqref{bcq1} holds for $h$ at $(\ox,\op,h(\ox,\op))$. Indeed, since the multifunction $G\colon p\mapsto\epi g(\cdot,p)$ is Lipschitz-like around $(\op,\ox,g(\ox,\op))$, there exist constants $\gg,\nu>0$ with
\begin{eqnarray}\label{e.9.3}
G(p_1)\cap\B_\nu\big((\ox,g(\ox,\op))\big)\subset G(p_2)+\gg d_1(p_1,p_2)\B_{X\times\R}\;\mbox{ for all }\;p_1,p_2\in\B_\nu(\op).
\end{eqnarray}
Define the multifunction $H\colon p\mapsto\epi h(\cdot,p)$ and select a neighborhood $W$ of $(\ox,h(\ox,\op))$ so that
$$
\Big(x,t-\frac{s}{2}\|x-\ox\|^2\Big)\in\B_\nu\big(\ox,g(\ox,\op)\big)\;\mbox{ for }\;(x,t)\in W.
$$
Taking any $p_1,p_2\in\B_\nu(\op)$ and $(x_1,t_1)\in H(p_1)\cap Z$, observe that $(x_1,t_1-\frac{s}{2}\|x_1-\ox\|^2)\in G(p_1)\cap\B_\nu((\ox,g(\ox,\op)))$, and thus we get by \eqref{e.9.3} that
\begin{eqnarray}\label{e.9.4}
\|x_2-x_1\|+\Big|r_2-t_1+\frac{s}{2}\|x_1-\ox\|^2\Big|\le\gg d_1(p_1,p_2)\;\mbox{ for some }\;(x_2,r_2)\in G(p_2).
\end{eqnarray}
Denoting $t_2:=r_2+\frac{s}{2}\|x_2-\ox\|^2$ yields $(x_2,t_2)\in H(p_2)$ and tells us together with \eqref{e.9.4} that
\begin{equation*}
\begin{array}{ll}
\|x_2-x_1\|+|t_2-t_1|&\disp\le\|x_2-x_1\|+\Big|r_2-t_1+\dfrac{s}{2}\|x_1-\ox\|^2\Big|+\dfrac{s}{2}\Big|\|x_2-\ox\|^2-\|x_1-\ox\|^2\Big|\\
&\disp\le\gg d_1(p_1,p_2)+\dfrac{s}{2}\Big|\|x_2-\ox\|-\|x_1-\ox\|\Big|\big(\|x_2-\ox\|+\|x_1-\ox\|\big)\\
&\disp\le\gg d_1(p_1,p_2)+\dfrac{s}{2}\|x_2-x_1\|\big(\|x_2-\ox\|+\|x_1-\ox\|\big)\\
&\disp\le\gg d_1(p_1,p_2)+\dfrac{s}{2}\gg d_1(p_1,p_2)2\nu=\gg(1+s\nu)d_1(p_1,p_2).
\end{array}
\end{equation*}
This readily brings us to the inclusion
\[
H(p_1)\cap W\subset H(p_2)+\gg(1+r\nu)d_1(p_1,p_2)\B_{X\times R}\;\mbox{ for all }\;p_1,p_2\in\B_\nu(\op),
\]
which verifies the claimed validity of BCQ \eqref{bcq1} for the function $h$ around $(\op,\ox,h(\ox,\op))$.

Employing now implication (i)$\Longrightarrow$(ii) in Lemma~\ref{fullmin} ensures the existence of a neighborhood $(U\times P\times V)$ of $(\ox,\op,\hat v)$ such that for all $(u,p,v)\in\gph\partial_x h\cap(U\times P\times V)$ we have
\begin{equation}\label{e.10}
\la z,w\ra\ge (s-r)\|w\|^2\quad\mbox{whenever}\quad z\in(\Hat D^*\partial h_p)(u,v)(w).
\end{equation}
To verify finally the second-order condition \eqref{4.14}, pick any $z\in(\Hat D^*\partial g_p)(u,v)(w)$ with $w\in X$ and $(u,p,v)\in W$; this implies that $(u,p,v+s(x-\ox))\in U\times P\times V$. It follows from \cite[Theorem~1.62]{M1} that $z+sw\in(\Hat D^*\partial h_p)(u,v+s(x-\ox))(w)$. Combining this with \eqref{e.10} yields $\la z+sw,w\ra\ge(s-r)\|w\|^2$, which clearly justifies \eqref{4.14} and hence completes the proof of the lemma.\endproof

Based on the above, we can derive now the first formula for calculating the threshold of prox-regularity that is valid in the case of Hilbert decision spaces. \vspace*{-0.1in}

\begin{Theorem}{\bf (threshold of prox-regularity in Hilbert spaces).}\label{thm3} Consider the potential $g$ of PVS \eqref{VS} in the setting of Lemma~{\rm\ref{lm0}} and define the number
\begin{equation}\label{e.8}
\tau:=\lim_{\eta\dn0}\left[\inf\Big\{\frac{\la z,w\ra}{\|w\|^2}\Big|\;z\in(\Hat D^*\partial g_p)(x,v)(w),\;(x,p,v)\in\gph\partial_x g\cap\B_\eta(\ox,\op,\hat v),\;w\ne 0\Big\}\right].
\end{equation}
Then the following assertions hold:
	
{\bf (i)} The number $\tau$ is finite.
	
{\bf (ii)}  The threshold of prox-regularity of $g$ at $(\ox,\op)$ for $\hat v$ is computed by $\mathcal{R}=\max\{0,-\tau\}$.
\end{Theorem}\vspace*{-0.05in}
{\bf Proof.} It follows from Lemma~\ref{lm0} and the definition of $\tau$ in \eqref{e.8} that $\tau\ge-r$ for any prox-parameter $r$ in \eqref{par-prox}. This $\tau>-\infty$ while being actually finite in the setting under consideration. It justifies (i) and show furthermore that $\mathcal{R}\ge\max\{0,-\tau\}$.

To verify (ii), choose any $s>-\tau$ and define $h(x,p):=g(x,p)+\frac{s}{2}\|x-\ox\|^2$ for $x\in X$ and $p\in\mathcal{P}$. We clearly have that $\partial h_p(x)= \partial g_p(x)+s(x-\ox)$ and
$$
\big(\Hat D^*\partial h_p\big)\big(x,v+s(x-\ox)\big)(w)=\big(\Hat D^*\partial g_p\big)(x,v)(w)+sw\;\mbox{ whenever }\;v\in\partial_x g(x,p).
$$
For any $\ve\in (0,s+\tau)$ we find some $\eta>0$ such that
\begin{equation}\label{zzx}
\la z,w\ra\ge (\tau-\ve) \|w\|^2\;\;\;\mbox{for all}\;\;\;z\in(\Hat D^*\partial g_p)(x,v)(w),\;(x,p,v)\in\gph\partial_x g\cap\B_\eta(\ox,\op,\hat v), w\in X.
\end{equation}
Fixing any pair $(x,v)\in\gph\partial g_p\cap\B_\eta(\ox,\hat v)$ and then picking $z\in(\Hat D^*\partial h_p)(x,v+s(x-\ox))(w)$ give us $z-sw\in(\Hat D^*\partial g_p)(x,v)(w)$. It follows from \eqref{zzx} that $\la z-sw,w\ra\ge(\tau-\ve)\|w\|^2$, which yields
\begin{equation}\label{e.10.1}
\la z,w\ra\ge(\tau+s-\ve)\|w\|^2\;\mbox{ for all }\;z\in\big(\Hat D^*\partial h_p\big)\big(x,v+s(x-\ox)\big)(w).
\end{equation}
Since $\tau+s-\ve>0$, combining \eqref{e.10.1} with implication (ii)$\Longrightarrow$(i) of Lemma~\ref{fullmin} tells us that USOGC \eqref{2.10} holds for $h$ at $(\ox,\op,\hat v)$ with modulus $\tau+s-\ve$. Hence there are neighborhoods $U_1$ of $\ox$, $P_1$ of $\op$, and $V_1$ of $\hat v$ on which we have the estimate
\begin{equation}\label{e.11}
h_p(x)\ge h_p(u)+\la v,x-u\ra+\dfrac{\tau+s-\ve}{2}\|x-u\|^2\;\mbox{ for all }\;u\in U_1,\;(u,p,v)\in\gph\partial h_{p}\cap(U_1\times P_1\times V_1).
\end{equation}

To complete the proof of this theorem, take any $(u,p,v)\in\gph\partial_x g\cap W$ with $W:=M(U_1\times P_1\times V_1)$ and $M(x,p,v):=(x,p,v-s(x-\ox))$. It follows from \eqref{e.11} that
\[\begin{array}{ll}
g_p(x)&\disp=h_p(x)-\dfrac{s}{2}\|x-\ox\|^2\\
&\disp\ge h_p(u)+\la v+s(u-\ox),x-u\ra+\frac{\tau+s-\ve}{2}\|x-u\|^2-\frac{s}{2}\|x-\ox\|^2\\
&\disp=g_p(u)+\la v,x-u\ra+\frac{\tau-\ve}{2}\|x-u\|^2\;\mbox{ for any }\;x\in U_1,
\end{array}\]
which justifies the prox-regularity of $g$ at $(\ox,\op)$ with modulus $\max\{0,\ve-\tau\}\ge\mathcal{R}$ for any $\ve\in (0,s+\tau)$. Combining it with the arguments in (i) tells us that $\mathcal{R}=\max\{0,-\tau\}$ and so ends the proof.
\endproof

The next theorem provides a {\em pointbased} formula for computing the threshold of prox-regularity of the potential $g$ in \eqref{VS} while assuming in addition the finite dimensionality of the decision space and the parameter independence of the potential.\vspace*{-0.1in}

\begin{Theorem}{\bf (threshold of prox-regularity in finite dimensions).}\label{thm4} Suppose that in the setting of Theorem~{\rm\ref{thm3}} we have $\dim X<\infty$ and $g=g(x)$ in \eqref{VS} with $\hat v\in\partial g(\ox)$. Then for
\begin{eqnarray}\label{e.12}
\tau_0:=\inf\Big\{\frac{\la z,w\ra}{\|w\|^2}\Big|\;z\in(D^*\partial g)(\ox,\hat v)(w),\;w\ne 0\Big\}
\end{eqnarray}
the following assertions are satisfied:
	
{\bf (i)} The number $\tau_0$ is well-defined.
	
{\bf (ii)} The threshold of prox-regularity of $g$ at $(\ox,\op)$ for $\hat v$ is computed by $\mathcal{R}=\max\{0,-\tau_0\}$.
\end{Theorem}\vspace*{-0.05in}
{\bf Proof.} To verify (i), take any modulus $r>0$ from definition \eqref{par-prox} of prox-regularity for $g$. Picking any $z\in D^*\partial g(\ox,\hat v)(w)$, deduce by passing to the limit in the necessary condition \eqref{4.14} of Lemma~\ref{lm0} that $\la z,w\ra\ge-r\|w\|^2$. This yields $\tau_0\ge-r$ and so justifies assertion (i) together with the lower estimate $\mathcal{R}\ge\max\{0,-\tau_0\}$ for the threshold of prox-regularity.

The proof of assertion (ii) follows the lines in the proof of Theorem~\ref{thm3} with appropriate modifications. Fix any $s>-\tau_0$ and define the function $h(x):=g(x)+\frac{s}{2}\|x-\ox\|^2$ as $x\in X$ for which $\partial h(x)=\partial g(x)+s(x-\ox)$. Then picking any $z\in(D^*\partial h)(\ox,\hat v)(w)$ gives us $z-sw\in(D^*\partial g)(\ox,\hat v)(w)$, and we conclude get from the definition of $\tau_0$ in \eqref{e.12} that
\begin{equation}\label{e.14.1}
\la z,w\ra\ge(\tau_0+s)\|w\|^2.
\end{equation}
Since $\tau_0+s>0$, it follows from \eqref{e.14.1} and the nonparametric pointbased version of Lemma~\ref{fullmin} from \cite[Theorems~3.2 and 3.6]{MN1} that USOGC \eqref{2.10} holds for $h$ at $(\ox,\hat v)$ with modulus $\tau_0+s-\ve>0$ if $\ve>0$ is sufficiently small. This allows us to find neighborhoods $U_1$ of $\ox$ and $V_1$ of $\hat v$ with
\begin{equation}\label{e.14.2}
h(x)\ge h(u)+\la v,x-u\ra+\frac{\tau_0+s-\ve}{2}\|x-u\|^2\;\mbox{ whenever }\;u\in U_1,\;(u,v)\in\gph\partial h\cap(U_1\times V_1).
\end{equation}
To finish the proof of the theorem, take any $(u,v)\in\gph\partial g\cap W$ with $W:=M(U_1\times V_1)$ for $M(x,v):=(x,v-s(x-\ox))$ and derive from
\eqref{e.14.2} that
\[\begin{array}{ll}
g(x)&\disp=h(x)-\frac{s}{2}\|x-\ox\|^2\\
&\disp\ge h(u)+\la v+s(u-\ox),x-u\ra+\frac{\tau_0+s-\ve}{2}\|x-u\|^2-\dfrac{s}{2}\|x-\ox\|^2\\
&\disp=g(u)+\la v,x-u\ra+\frac{\tau_0-\ve}{2}\|x-u\|^2,
\end{array}\]
where the inequality comes from \eqref{e.14.2} with $v+s(u-\ox)\in\partial h(u)$. This verifies the prox-regular at $\ox$ with modulus $\max\{0,\ve-\tau_0\}$. Since $\ve>0$ was chosen to be arbitrarily small, we get the upper estimate $\mathcal{R}\le\max\{0,-\tau_0\}$ and hence completes the proof of (ii) and of the whole theorem.\endproof\vspace*{-0.15in}

\begin{Remark}{\bf (on hypomonotonicity).}\label{rm5} {\rm It is worth mentioning that the prox-regularity of any extended-real-valued l.s.c.\ function is equivalent and quantitatively related to {\em hypomonotonicity} property of its subdifferential; see \cite[Theorem~13.36, Example~12.28]{rw} and the discussions therein. The crucial roles of hypomonotonicity and threshold $\mathcal{R}$ of prox-regularity have been recognized not only in variational theory but also in various optimization algorithms as, e.g., the proximal point method in Pennanen \cite{P}. It seems that our paper is the first in the literature with the {\em exact computations} of the prox-regularity threshold $\mathcal{R}$, which is a crucial quantitative characteristics.}
\end{Remark}\vspace*{-0.25in}

\section{H\"olderian and Lipschitzian Properties of Proximal Mappings}
\setcounter{equation}{0}

From one hand, this section can be considered as a preliminary step to derive the main results on H\"olderian and Lipschitzian full stability of PVS \eqref{VS} without any smoothness assumptions on the base mapping $f$. On the other hand, the results obtained here on H\"olderian and Lipschitzian properties of the mapping $P_\lm\colon X\times\mathcal{P}\tto X$ given by
\begin{eqnarray}\label{e.1}
P_\lm(v,p):=\big\{x\in X|\;v\in x+\lm\partial_x g(x,p)\big\},\quad\lm>0,
\end{eqnarray}
certainly are of their own interest and importance; see below. If $g$ is convex with respect of $x$, \eqref{e.1} agrees with the {\em proximal mapping} of $g$; we'll keep this name in the general case under consideration. If furthermore $g(x,p)=\delta_{C(p)}(x)$ is the indicator function of the parameter-dependent set $C(p)$, \eqref{e.1} reduces to the {\em projection mapping} associate with $C(p)$.

Observe that $P_\lm$ is the solution map for the following PVS of type \eqref{VS} with the smooth parameter-independent base mapping $f(x,p,q)=x$:
\begin{equation}\label{e.1.0}
v\in x+\partial_x(\lm g)(x,p),\quad\lm>0.
\end{equation}
Although general characterizations of full stability for the smooth-base PVS \eqref{VS} and their specification in \eqref{e.1.0} are given in Section~4, we need in what follows a certain {\em quantitative} version, which relates the threshold ${\cal R}$ of prox-regularity of the potential $g$ and the modulus $\sigma>0$ of the {\em local strong monotonicity} of the base $f$ in \eqref{VS} with respect to $x$ defined as
\begin{eqnarray}\label{smon}
\la f(x_1,p,q)-f(x_2,p,q),x_1-x_2\ra\ge\sigma\|x_1-x_2\|^2\;\mbox{ for any }\;x_1,x_2\in U,\;(p,q)\in P\times Q
\end{eqnarray}
on some neighborhood $U\times P\times Q$ of the reference triple $(\ox,\op,\oq)$. This is done in the next statement, which is established with the help of Lemma~\ref{lm0} under the additional partial smoothness assumption (A1) on the base $f$ imposed in Section~4.\vspace*{-0.1in}

\begin{Proposition}{\bf (quantitative condition for H\"olderian full stability of smooth PVS).}\label{prop0} Let $\ox\in S(\ov,\op,\oq)$ in the setting of Theorem~{\rm\ref{Holderpq}} with some prox-parameter $r>0$ in \eqref{par-prox} and the threshold ${\cal R}$ of prox-regularity of the potential $g$ in {\rm(A2)}. Suppose in addition that the base mapping $f$ is locally strongly monotone with respect to the decision variable $x$ around $(\ox,\op,\oq)$, and the strong monotonicity modulus satisfies $\sigma>{\cal R}$. Then the second-order condition \eqref{4.8a} in Theorem~{\rm\ref{Holderpq}(ii)} holds with $\kk_0:=\sigma-r$, and hence $\ox$ is a H\"olderian fully stable solution of PVS \eqref{VS}.
\end{Proposition}\vspace*{-0.07in}
{\bf Proof.} Since $\sigma>R$, we can take $r\in(\mathcal{R},\sigma)$ for the prox-parameter $r$ in \eqref{par-prox}. It is easy to see that the strong monotonicity of the base $f$ in \eqref{smon} and its smoothness in $x$ ensure that
\begin{eqnarray}\label{4.15}
\la\nabla_x f(\ox,\op,\oq)w,w\ra\ge\sigma\|w\|^2\;\mbox{ for all}\;w\in X.
\end{eqnarray}
Applying now Lemma~\ref{lm0}, we find a neighborhood $W$ of $(\ox,\op,\hat v)$ on which condition \eqref{4.14} is satisfied. It implies together with \eqref{4.15} that the estimate
$$
\la\nabla_x f(\ox,\op,\oq)w,w\ra+\la z,w\ra\ge\sigma\|w\|^2-r\|w\|^2=\kk_0\|w\|^2
$$
holds for all $(w,z)\in X\times X$ with $z\in(\Hat D^*\partial g_p)(u,v)(w)$, where $\kk_0=\sigma-r>0$. This gives us condition {\eqref{4.8a} of
Theorem~\ref{Holderpq}(ii) and hence verifies the claimed H\"olderian full stability.\endproof

The obtained proposition allows us to derive the following lemma, which is instrumental to establish the main results of this and next sections without any smoothness assumptions. For the convenience in further references, we formulate this result in the general framework of PVS \eqref{VS} while it primarily concerns the proximal mapping \eqref{e.1} depending only on the potential $g$ of \eqref{VS}.\vspace*{-0.1in}

\begin{Lemma}{\bf (H\"olderian localization of the proximal mapping).}\label{prop1} Let $\ox\in S(\ov,\op,\oq)$ in \eqref{ss} under the validity of assumptions {\rm(A2)} and {\rm(A3)} with some prox-parameter $r>0$ in \eqref{par-prox} and the threshold ${\cal R}$ of prox-regularity of $g$ in {\rm(A2)}. Then for any $\lm\in(0,\mathcal{R}^{-1})$ the proximal mapping $P_\lm$ has a single-valued localization $\pi_\lm$ relative to a neighborhood $W\times P\times U$ of
$(\lm\hat v+\ox,\op,\ox)$ such that for all $(v_1,p_1),(v_2,p_2)\in W\times P$ we have
\begin{equation}\label{e.1.1}
\|(v_1-v_2)-2\kk_0[\pi_\lm(v_1,p_1)-\pi_\lm(v_2,p_2)]\|\le\|v_1-v_2\|+\ell_0 d_1(p_1,p_2)^\frac{1}{2}
\end{equation}
with $\kk_0:=1-\lm r$ and some positive constant $\ell_0$.
\end{Lemma}\vspace*{-0.05in}
{\bf Proof.} Fix $\lm\in(0,\mathcal{R}^{-1})$ and let $\tilde{v}:=\lm\hat v+\ox$. It follows from (A2) with a given prox-parameter $r$ that there exists a neighborhood $U\times P\times V$ of $(\ox,\op,\hat v)$ on which \eqref{par-prox} are satisfied. By the choice of $\lm$ we can always suppose that $0<\lm<r^{-1}<{\cal R}^{-1}$. Considering now the proximal PVS \eqref{e.1.0}, observe that assumptions (A2) and (A3) hold for this system, where the imposed parametric continuous prox-regularity of $g$ at $(\ox,\op)$ for $\hat v$ yields this property for $\lm g$ at $(\ox,\op)$ for $\tilde{v}$ with the prox-parameter $\tilde r:=\lm r<1$. Thus $\tilde{\cal R}<1$ for the threshold of prox-regularity of the potential $\lm g$ of \eqref{e.1.0} at $(\ox,\op)$ for $\tilde{v}$. Furthermore, the base mapping $f(x,p,q)=x$ in \eqref{e.1.0} is strongly monotone on $X$ with modulus $\sigma=1$, and so $\sigma>\tilde{\cal R}$. Applying Proposition~\ref{prop0} ends the proof of the lemma.\endproof

The following theorem shows that the mapping $P_\lm$ has a single-valued localization satisfying a H\"olderian estimate expressed in terms of the Pompeiu-Hausdorff distance between parametric epigraphical sets for $g$. It has various important subsequences one of which is discussed in Example~\ref{ex0} for the case where $g$ is the indicator function of a parameter-dependent convex set.\vspace*{-0.1in}

\begin{Theorem}{\bf (H\"olderian estimate for the proximal mapping via the Pompeiu-Hausdorff distance).}\label{thm1} In the setting of Lemma~{\rm\ref{prop1}} we have that the proximal mapping $P_\lm$ has a single-valued localization $\pi_\lm$ relative to a neighborhood $W\times P\times U$ of $(\lm\hat v+\ox,\op,\ox)$ such that
\begin{equation}\label{e.15}
\|(v_1-v_2)-2(1-\lm r)[\pi_\lm(v_1,p_1)-\pi_\lm(v_2,p_2)]\|\le\|v_1-v_2\|+\ell_1\th+\ell_2\th^{\frac{1}{2}}
\end{equation}
for all $(v_1,p_1),(v_2,p_2)\in W\times P$ with some explicitly calculated constants $\ell_1,\ell_2>0$, where
\begin{eqnarray}\label{th}
\th:={\rm haus}\big[\epi g(\cdot,p_1)\cap(U\times\R),\epi g(\cdot,p_2)\cap(U\times\R)\big].
\end{eqnarray}
\end{Theorem}\vspace*{-0.05in}
{\bf Proof.} Employing Lemma~\ref{prop1}, for any $\lm\in(0,\mathcal{R}^{-1})$ we get the single-valued localization $\pi_\lm$ of $P_\lm$ satisfying condition \eqref{e.1.1} in the neighborhood $W\times P\times U$ with $\kk_0=1-\lm r$ and some constant $\ell_0>0$. Taking any $v_1,v_2\in U,$ $p_1,p_2\in P$, $x_i\in P_\lm(v_i,p_i)\cap U$ for $i=1,2$, and $x_3\in P_\lm(v_1,p_2)\cap U$, assume without loss of generality that $\hat{v}=0$ with $\ov=\ox$, $U:=\B_{\rho}(\ox)$ for some constant $\rho>0$, and that $V:=\B_\nu(\hat{v})$ with $0<\rho\max\{1,2\lm^{-1}\}<\nu$. For any $i\in\{1,2\}$ and $j\in\{1,2,3\}$ we get
$$
\|\lm^{-1}(v_i-x_j)\|\le\lm^{-1}(\|v_i-\ox\|+\|\ox-x_j\|)\le 2\lm^{-1}\rho<\nu,
$$
and hence $\lm^{-1}(v_i-x_j)\in V$ for such $i,j$. It easily follows from estimate \eqref{e.1.1} that
\begin{eqnarray}\label{e.16}
\|v_1-v_2-2(1-\lambda r)(x_3-x_2)\|\le\|v_1-v_2\|.
\end{eqnarray}
Let us now show that there exist positive constants $\ell_1$ and $\ell_2$ such that
\begin{eqnarray}\label{e.17}
2(1-\lambda r)\|x_1-x_3\|\le\ell_1\th+\ell_2\th^\frac{1}{2}.
\end{eqnarray}
To proceed, let $\alpha_1:=g(x_1,p_1)$ and $\alpha_3:=g(x_3,p_2)$ for which $(x_1,\alpha_1)\in\epi g(\cdot,p_1)\cap(U\times\R)$ and $(x_3,\alpha_3)\in\epi g(\cdot,p_2)\cap(U\times\R)$. The lower semicontinuity of $g$ yields the closedness of the sets $\epi g(\cdot,p_i)\cap(U\times\R)$, $i=1,2$, with the closed ball $U$ defined above. This allows us to find $(y_1,s_1)\in\epi g(\cdot,p_2)\cap(U\times\R)$ and $(y_3,s_3)\in\epi g(\cdot,p_1)\cap(U\times\R)$ such that
\[(y_1,s_1)=(x_1,\alpha_1)+(z_1,r_1),\;(y_3,s_3)=(x_3,\alpha_3)+(z_3,r_3)\;\mbox{ with }\;\|(z_i,r_i)\|\le\th,\;i=1,3,
\]
by the definition of $\th$ in \eqref{th}. Since $\lm^{-1}(v_1-x_i)\in V$ as $i=1,3$, applying \eqref{par-prox} yields
\begin{eqnarray*}
g(y_3,p_1)\ge g(x_1,p_1)+\la\lm^{-1}(v_1-x_1),y_3-x_1\ra-\frac{r}{2}\|y_3-x_1\|^2,\\
g(y_1,p_2)\ge g(x_3,p_2)+\la\lm^{-1}(v_1-x_3),y_1-x_3\ra-\frac{r}{2}\|y_1-x_3\|^2.
\end{eqnarray*}
Adding these inequalities gives us the relationships
\begin{equation*}
\begin{array}{ll}
&\disp g(y_3,p_1)+g(y_1,p_2)-g(x_1,p_1)-g(x_3,p_2)\\
&\disp\ge\lm^{-1}\la v_1-x_1,y_3-x_1\ra+\lm^{-1}\la v_1-x_3,y_1-x_3\ra-\frac{r}{2}\|y_3-x_1\|^2-\frac{r}{2}\|y_1-x_3\|^2\\
&\disp=\lm^{-1}\la v_1-x_1,x_3-x_1\ra+\lm^{-1}\la v_1-x_3,x_1-x_3\ra+\lm^{-1}\la v_1-x_1,z_3\ra+\lm^{-1}\la v_1-x_3,z_1\ra\\
&\disp\qquad\qquad-\frac{r}{2}\left(\|z_1\|^2+\|z_3\|^2+2\la z_1-z_3,x_1-x_3\ra+2\|x_1-x_3\|^2\right)\\
&\disp=(\lm^{-1}-r)\|x_1-x_3\|^2+\lm^{-1}\left(\la v_1-x_1,z_3\ra+\la v_1-x_3,z_1\ra\right)\\
&\disp\qquad\qquad-\frac{r}{2}\left(\|z_1\|^2+\|z_3\|^2\right)-r\la z_1-z_3,x_1-x_3\ra.
\end{array}
\end{equation*}
By using $\|(z_i,r_i)\|\le\th$ for $i=1,3$ we have the estimates
\begin{eqnarray*}
\begin{array}{ll}
g(y_3,p_1)+g(y_1,p_2)-g(x_1,p_1)-g(x_3,p_2)\le s_3+s_1-\alpha_1-\alpha_3=r_1+r_3\le 2\th;\\
\lm^{-1}\left(\la v_1-x_1,z_3\ra+\la v_1-x_3,z_1\ra\right)=\lm^{-1}\left(\la v_1-x_1,z_1+z_3\ra+\la x_1-x_3,z_3\ra\right)\\
\ge-\lm^{-1}\|v_1-x_1\|\cdot\|z_1+z_3\|-\lm^{-1}\th\|x_1-x_3\|\ge-4\lm^{-1}\rho\th-\lm^{-1}\th\|x_1-x_3\|;\\
-\dfrac{r}{2}\left(\|z_1\|^2+\|z_3\|^2\right)\ge-r\th^2\ge-\lm^{-1}\th^2;\quad\mbox{and}\\
-r\la z_1-z_3,x_1-x_3\ra\ge-2r\th\|x_1-x_3\|\ge-2\lm^{-1}\th\|x_1-x_3\|.
\end{array}
\end{eqnarray*}
As a result, it follows from the above that
\begin{eqnarray*}
0\ge(\lm^{-1}-r)\|x_1-x_3\|^2-3\lm^{-1}\th\|x_1-x_3\|-4\lm^{-1}\rho\th-\lm^{-1}\th^2-2\th,
\end{eqnarray*}
which can be equivalently rewritten as
\[
0\ge(1-\lm r)\|x_1-x_3\|^2-3\th\|x_1-x_3\|-4\rho\th-2\lm\th-\th^2.
\]
The right-hand side of this inequality is a quadratic expression of $\|x_1-x_3\|$ with the discriminant
\begin{eqnarray*}
\Delta:=9\th^2-4(1-\lm r)(-4\rho\th-2\lm\th-\th^2)=(13-4\lm r)\th^2+8(2\rho+\lm)(1-\lm r)\th>0.
\end{eqnarray*}
By solving the quadratic equation above with taking into account that $1-\lm r>0$, we arrive at
\begin{equation*}
\|x_1-x_3\|\le\frac{1}{2(1-\lm r)}\left(3\th+\big\{(13-4\lm r)\th^2+8(2\rho+\lm)(1-\lm r)\th\big\}^\frac{1}{2}\right),
\end{equation*}
which implies in turn the inequalities
\begin{equation*}
2\kk_0\|x_1-x_3\|\le 3\th+\big\{(9+4\kk_0)\th^2+8(2\rho+\lm)\kk_0\th\big\}^\frac{1}{2}\le\left(3+\sqrt{9+4\kk_0}\right)\th+2\sqrt{2(2\rho+\lm)\kk_0}\th^\frac{1}{2}
\end{equation*}
with $\kk_0=1-\lm r$. This verifies \eqref{e.17} with the {\em explicitly calculated} constants
\begin{eqnarray}\label{ell}
\ell_1:=3+\sqrt{9+4\kk_0}\quad\mbox{and}\quad\ell_2:=2\sqrt{2(2\rho+\lm)\kk_0}.
\end{eqnarray}
Combining finally \eqref{e.16} and \eqref{e.17} gives us the estimates
\[
\|(v_1-v_2)-2\kk_0(x_1-x_2)\|\le\|(v_1-v_2)-2\kk_0(x_3-x_2)\|+2\kk_0\|x_1-x_3\|\le\|v_1-v_2\|+\ell_1\th+\ell_2\th^{\frac{1}{2}}
\]
with the constants $\kk_0,\ell_1,\ell_2$ calculated above. It justifies \eqref{e.15} and hence completes the proof.
\endproof

Let us now specify the results of Theorem~\ref{thm1} for the case where the potential $g$ in \eqref{VS} is the indicator function of a parameter-dependent (may not be convex) set. This setting was considered by Robinson \cite{R05} when both decision and parameter spaces are finite-dimensional.\vspace*{-0.1in}

\begin{Example}{\bf(H\"olderian properties of projections on prox-regular sets).}\label{ex0}
{\rm Let $\partial_x g(x,p)=N_{C(p)}(x)$ in \eqref{VS} via the limiting normal cone to the closed-valued mapping $C\colon\mathcal{P}\tto X$. In this case the subdifferential continuity in (A2) holds automatically while (A3) reduces to the Lipschitz-like property of $C$ around the reference point $\op$. Note that the variational system \eqref{e.1.0} is written now as $v\in x+N_{C(p)}(x)$ and the proximal mapping \eqref{e.1} reads as
$$
P_\lm(v,p)=\big\{x\in X\big|\;x\in\big(N_{C(p)}+I\big)^{-1}(v)\big\}=\big(N_{C(p)}+I\big)^{-1}(v),
$$
which is the projector on $C(p)$. It is clear that the $\lm$-parameter is not needed in this setting. Observe also that in this case the quantity $\th$ from \eqref{th} reduces to $\th={\rm haus}[C(p_1)\cap U,C(p_2)\cap U]$ and estimate \eqref{e.15} of Theorem~\ref{thm1} implies, with $v_1=v_2=:v$, that
\begin{eqnarray}\label{hol-pr}
\|\pi_\lm(v,p_1)-\pi_\lm(v,p_2)\|\le\frac{\ell_1}{2(1-r)}\th+\frac{\ell_2}{2(1-r)}\th^{\frac{1}{2}},
\end{eqnarray}
where the positive constants $\ell_1$ and $\ell_2$ are explicitly calculated in the proof of Theorem~\ref{thm1}.

Note that Robinson's result in \cite[Theorem~3]{R05} provides a quantitative estimate for a single-valued projection localization of a prox-regular moving set $C(p)\subset\R^n$ that merely continuously depends on its parameter vs. the Lipschitz-like property of $C(p)$ corresponding to BCQ \eqref{bcq1} in (A3). The aforementioned estimate in \cite{R05} is similar albeit different from \eqref{hol-pr} and is derived differently from the proof of Theorem~\ref{thm1} given above.

To conclude the discussion in this example, we note further that when the set-valued mapping $C$ has closed convex values, our  Theorem~6.2 clearly implies the H\"olderian stability of the projection mapping obtained by Yen \cite[Lemma~1.1]{Y2} in which his Lipschitz-like assumption on the mapping $C$ is exactly our BCQ for $g(x,p)=\delta_{C(p)}(x)$; see \eqref{haus} and further analysis in our Remark~7.4.}
\end{Example}\vspace*{-0.1in}

The final result of this section establishes an effective {\em characterization} of the existence of a Lipschitzian single-valued localization of the proximal mapping \eqref{e.1} via the Lipschitz-like property of the graphical subdifferential mapping \eqref{K} imposed and employed above.\vspace*{-0.1in}

\begin{Theorem}{\bf (Lipschitzian localization of the proximal mapping).}\label{prop2} The following assertions are equivalent in the setting of Lemma~{\rm\ref{prop1}}:

{\bf (i)} Given any $\lm\in(0,\mathcal{R}^{-1})$, the proximal mapping \eqref{e.1} has a single-valued localization $\pi_\lm$ relative to a neighborhood $W\times P\times U$ of $(\lm \hat v+\ox,\op,\ox)$ so that the Lipschitzian estimate
\begin{eqnarray}\label{e.7.1}
\|(v_1-v_2)-2\kk_0[\pi_\lm(v_1,p_1)-\pi_\lm(v_2,p_2)]\|\le\|v_1-v_2\|+\ell_0 d_1(p_1,p_2)
\end{eqnarray}
holds for all $(v_1,p_1),(v_2,p_2)\in W\times P$ with some positive constants $\kk_0$ and $\ell_0$.

{\bf (ii)} The graphical subdifferential mapping \eqref{K} is Lipschitz-like around $(\op,\ox,\hat v)$.
\end{Theorem}\vspace*{-0.07in}
{\bf Proof.} It follows from Lemma~\ref{prop1} that for any $\lm\in(0,\mathcal{R}^{-1})$ the point $\ox$ is a H\"olderian fully stable solution to the proximal PVS \eqref{e.1.0} corresponding to the parameter pair $(\ov,\op)$. Employing the characterization in Theorem~\ref{Holderpq} ensures the existence of positive numbers $\eta,\kk_0$ such that whenever $(u,p,v)\in\gph\partial_x(\lm g)\cap\B_\eta(\ox,\op,\hat v)$ with $\hat v=\ov-\ox$ we have
\begin{eqnarray*}
\la w,w\ra+\la z,w\ra\ge\kk_0\|w\|^2\quad\mbox{for all}\quad z\in\big(\Hat D^*\partial(\lm g_p)\big)(u,v)(w),\;w\in X.
\end{eqnarray*}
Assuming now the validity of (ii) gives us by Proposition~\ref{Lippq} that $\ox$ is Lipschitzian fully stable for PVS \eqref{e.1.0}, and hence \eqref{e.7.1} holds
by Definition~\ref{fs}(ii). This verifies implication (ii)$\Longrightarrow$(i).

To verify the converse implication (i)$\Longrightarrow$(ii), observe that the assumed estimate \eqref{e.7.1} in (i) exactly constitutes the Lipschitzian full stability of $\ox$ in the proximal PVS \eqref{e.1.0}, and so the Lipschitz-like property in (ii) immediately follows from Proposition~\ref{Lippq} applied to this system.\endproof\vspace*{-0.2in}

\section{Full Stability of Nonsmooth PVS under Strong Monotonicity}
\setcounter{equation}{0}

In this section we study H\"olderian and Lipschitzian full stability of the general PVS \eqref{VS} without any smoothness assumptions on the base mapping $f$ while with imposing instead of the strong monotonicity property \eqref{smon} on this mapping with respect to the state variable $x$. Our goal is to find a relationship between the strong monotonicity modulus $\sigma>0$ and the threshold of prox-regularity of the potential function $g$ in \eqref{VS} ensuring the validity of both H\"olderian and Lipschitzian full stability of such systems under the standing assumptions (A1)--(A4) on the given data. The H\"olderian result of this type was obtained in Proposition~\ref{prop0} under the smoothness assumption (A1) on the base, which is significantly employed in the proof. In what follows we develop a different approach to full stability, which implements the result of Section~6 and allows us to completely drop (A1).

Our first result here concerns the H\"olderian version of full stability for \eqref{VS}.\vspace*{-0.1in}

\begin{Theorem}{\bf (H\"olderian full stability under strong monotonicity).}\label{thm2} Let $\ox\in S(\ov,\op,\oq)$ in \eqref{ss} under the standing assumptions {\rm(A1)}--{\rm(A3)} with the threshold ${\cal R}$ of prox-regularity of $g$ in {\rm(A2)}. Assume in addition that the base mapping $f$ is locally strongly monotone in $x$ on some neighborhood $U\times P\times Q$ of $(\ox,\op,\oq)$ with modulus $\sigma>{\cal R}$ in \eqref{smon}. Then $\ox$ is a H\"olderian fully stable solution to \eqref{VS} corresponding to the parameter triple $(\ov,\op,\oq)$.
\end{Theorem}\vspace*{-0.05in}
{\bf Proof.} Suppose without loss of generality that $r\in({\cal R},\sigma)$ for the prox-parameter $r$ in \eqref{par-prox} with the same neighborhoods $U$ and $P$ as in \eqref{smon}. Lemma~\ref{prop1} allows us to find a single-valued localization $\pi_\lm$ of the proximal mapping $P_\lm$ at $(\lm\hat v+\ox,\op)$ for $\ox$ for any $\lm\in(0,r^{-1})$. Thus there is a neighborhood $U_1\times P_1\times V_1\subset U\times P\times X$ of $(\ox,\op,\tilde v)$ with $\tilde v:=\lm \hat v+\ox$ such that $\gph\pi_\lm=\gph P_\lm\cap (V_1\times P_1\times U_1)$ and that estimate \eqref{e.1.1} holds for all $(v_1,p_1),(v_2,p_2)\in V_1\times P_1$ with $\kk_0=1-r\lm$ and some constant $\ell_0>0$. This implies in turn the inequality
\begin{equation}\label{e.5.1}
\|\pi_\lm(v_1,p_1)-\pi_\lm(v_2,p_2)\|\le\frac{1}{\kk_0}\|v_1-v_2\|+\frac{\ell_0}{2\kk_0}d_1(p_1,p_2)^{\frac{1}{2}}.
\end{equation}
Observe furthermore that the feasibility condition $x\in S(v,p,q)$ can be equivalently written as
\[
x+\lm v-\lm f(x,p,q)\in\lm\partial_x g(x,p)+x,\;\mbox{ i.e. }\;x\in P_\lm(x+\lm v-\lm f(x,p,q),p).
\]
It follows from the Lipschitz continuity of $f$ in (A1) that there is a neighborhood $U_2\times V_2\times P_2\times Q_2\subset U_1\times V_1\times P_1\times Q$ of $(\ox,\ov,\op,\oq)$ such that $(x+\lm v-\lm f(x,p,q),p)\in V_1\times P_1$ for all $(x,v,p,q)\in U_2\times V_2\times P_2\times Q_2$. Thus $x\in S(v,p,q)$ if and only if
\begin{eqnarray}\label{H}
x=\pi_\lm(x+\lm v-\lm f(x,p,q),p)=:H_{(v,p,q)}(x)\;\mbox{ for all }\;(x,v,p,q)\in U_2\times V_2\times P_2\times Q_2.
\end{eqnarray}
Fixing $(v,p,q)\in V_2\times P_2\times Q_2$, we claim that $H_{(v,p,q)}$ satisfies the {\em contraction condition} for some $\lm\in(0,r^{-1})$. Indeed, for any $x_1,x_2\in U_2$ it follows from the H\"older continuity \eqref{e.5.1} of the localization $\pi_\lm$ combined with \eqref{smon} and \eqref{lipf} that
\[\begin{array}{ll}
&\disp\|H_{(v,p,q)}(x_1)-H_{(v,p,q)}(x_2)\|^2=\|\pi_\lm(x_1+\lm v-\lm f(x_1,p,q),p)-\pi_\lm(x_2+\lm v-\lm f(x_2,p,q),p)\|^2\\
&\disp\le\frac{1}{(1-r\lm)^2}\|(x_1-x_2)-\lm(f(x_1,p,q)-f(x_2,p,q))\|^2\\&\disp=\frac{1}{(1-r\lm)^2}\left[\|x_1-x_2\|^2-2\lm\la f(x_1,p,q)-f(x_2,p,q),x_1-x_2\ra+\lm^2\|f(x_1,p,q)-f(x_2,p,q)\|^2\right]\\&\disp\le\frac{1}{(1-r\lm)^2}\left(1-2\lm\sigma+\lm^2L^2\right)\|x_1-x_2\|^2
=\left[1-\frac{\lm\big(2(\sigma-r)-\lm(L^2-r^2)\big)}{(1-r\lm)^2}\right]\|x_1-x_2\|^2,
\end{array}\]
where $L$ is a Lipschitz constant from (A1). Since $\sigma>r$, it is possible to choose $\lm\in(0,r^{-1})$ sufficiently small so that $2(\sigma-r)>\lm(L^2-r^2)$. By the above inequalities we have
\begin{eqnarray}\label{e.6}
\|H_{(v,p,q)}(x_1)-H_{(v,p,q)}(x_2)\|\le\al\|x_1-x_2\|\;\mbox{ with }\;\al:=\frac{\left(1-2\lm\sigma+\lm^2L^2\right)^\frac{1}{2}}{1-r\lm}<1.
\end{eqnarray}
Furthermore, the continuity of the mappings $\pi_\lm$ and $f$ allows us to find constants $\delta,\eta>0$ such that $\B_\delta(\ox)\times\B_\eta(\ov)\times\B_\eta(\op)\times\B_\eta(\oq)\subset U_2\times V_2\times P_2\times Q_2$ and that
\[
\|H_{(v,p,q)}(\ox)-\ox\|=\|\pi_\lm(\ox+\lm v-\lm f(\ox,p,q),p)-\pi_\lm(\ox +\lm v-\lm f(\ox,\op,\oq),\op)\|\le\delta(1-\al)
\]
for all $(v,p,q)\in\B_\eta(\ov)\times\B_\eta(\op)\times\B_\eta(\oq)$. This together with \eqref{e.6} tells us that
\[
\|H_{(v,p,q)}(x)-\ox\|\le\|H_{(v,p,q)}(x)-H_{(v,p,q)}(\ox)\|+\|H_{(v,p,q)}(\ox)-\ox\|<\al\|x-\ox\|+\delta(1-\alpha)\le\delta
\]}
for all $x\in\B_\delta(\ox)$. Then the classical Banach contraction theorem ensures that the mapping $H_{(v,p,q)}$ from \eqref{H} has a unique {\em fixed point} in $\B_\delta(\ox)$, denoted by $\vt(v,p,q)$, for which we have
\[
\gph\vt=\gph S\cap\big(\B_\eta(\ov)\times\B_\eta(\op)\times\B_\eta(\oq)\times\B_\delta(\ox)\big).
\]

To verify the claimed H\"olderian full stability, it remains to check the validity of \eqref{4.7}. We proceed by taking any $(v_1,p_1,q_1),(v_2,p_2,q_2)\in \B_\eta(\ov,\op,\oq)$ with $x_1:=\vt (v_1,p_1,q_1)$, $x_2:=\vt(v_2,p_2,q_2)$, and $x_3:=\vt(v_1,p_2,q_2)$ and observing that
$$
x_i=\pi_\lm(x_i+\lm v_i-\lm f(x_i,p_i,q_i),p_i),\;i=1,2,\;\mbox{ and }\;x_3=\pi_\lm(x_3+\lm v_1-\lm f(x_3,p_2,q_2),p_2).
$$
Employing \eqref{e.1.1} for the pairs $(x_2+\lm v_2-\lm f(x_2,p_2,q_2),p_2)$ and $(x_3+\lm v_1-\lm f(x_3,p_2,q_2),p_2)$ yields
\[\begin{array}{ll}
\|x_3-x_2+\lm (v_1-v_2)-\lm (f(x_3,p_2,q_2)-f(x_2,p_2,q_2))-2\kk_{0}(x_3-x_2)\|\\
\disp \le \|x_3-x_2+\lm (v_1-v_2)-\lm (f(x_3,p_2,q_2)-f(x_2,p_2,q_2))\|,
\end{array}\]
which implies in turn the following relationships:
\[\begin{array}{ll}
0\disp\le\|x_3-x_2+\lm (v_1-v_2)-\lm(f(x_3,p_2,q_2)-f(x_2,p_2,q_2))\|^2\\
\qquad\qquad\qquad-\|x_3-x_2+\lm (v_1-v_2)-\lm(f(x_3,p_2,q_2)-f(x_2,p_2,q_2))-2\kk_0 (x_3-x_2)\|^2\\
\disp=4\kk_0\la x_3-x_2+\lm(v_1-v_2)-\lm(f(x_3,p_2,q_2)-f(x_2,p_2,q_2)),x_3-x_2\ra-4\kk_0^2\|x_3-x_2\|^2\\
=4\kk_0(1-\kk_0)\|x_3-x_2\|^2+4\kk_0\lm\la v_1-v_2, x_3-x_2\ra-4\kk_0\lm\la f(x_3,p_2,q_2)-f(x_2,p_2,q_2),x_3-x_2\ra\\
\le 4\kk_0(1-\kk_0)\|x_3-x_2\|^2+4\kk_0\lm\la v_1-v_2, x_3-x_2\ra-4\kk_0\lm\sigma\|x_3-x_2\|^2\\
=4\kk_0(1-\kk_0-\lm\sigma)\|x_3-x_2\|^2+4\kk_0\lm\la v_1-v_2,x_3-x_2\ra,
\end{array}\]
where the last inequality follows from \eqref{smon}. Hence we have
\[
0\le\la v_1-v_2,x_3-x_2\ra-\kk\|x_3-x_2\|^2\;\mbox{ with }\;\kk:=\lm^{-1}(\kk_0+\lm\sigma-1)=\lm^{-1}(\lm\sigma-\lm r)=\sigma-r>0.
\]
This readily implies the relationships
\[\begin{array}{ll}
\|v_1-v_2-2\kk(x_3-x_2)\|^2=\|v_1-v_2\|^2-4\kk\big(\la v_1-v_2,x_3-x_2\ra-\kk\|x_3-x_2\|^2\big)\le\|v_1-v_2\|^2
\end{array}\]
and therefore the following estimate:
\begin{eqnarray}\label{e.7}
\|v_1-v_2-2\kk(x_3-x_2)\|\le\|v_1-v_2\|.
\end{eqnarray}
To proceed further, deduce from by the H\"older continuity of $\pi_\lm$ in \eqref{e.5.1} that
\begin{equation}\label{e.7.0a}
\begin{array}{ll}
\|x_1-x_3\|&\disp=\|H_{(v_1,p_1,q_1)}(x_1)-H_{(v_1,p_2,q_2)}(x_3)\|\\
&\disp=\|\pi_\lm(x_1+\lm v_1-\lm f(x_1,p_1,q_1),p_1)-\pi_\lm(x_3+\lm v_1-\lm f(x_3,p_2,q_2),p_2)\|\\
&\disp\le\frac{1}{1-r\lm}\|x_1-x_3-\lm (f(x_1,p_1,q_1)-f(x_3,p_2,q_2))\|+\ell d_1(p_1,p_2)^{\frac{1}{2}}\\
\end{array}
\end{equation}
with $\ell:=\frac{\ell_0}{2\kk_0}$. It follows furthermore that
\[\begin{array}{ll}
&\|x_1-x_3-\lm(f(x_1,p_1,q_1)-f(x_3,p_2,q_2))\|\\
&\disp\le\|x_1-x_3-\lm(f(x_1,p_1,q_1)-f(x_3,p_1,q_1))\|+\lm\|f(x_3,p_2,q_2)-f(x_3,p_1,q_1)\|\\
&\disp=\left[\|x_1-x_3\|^2-2\lm\la f(x_1,p_1,q_1)-f(x_3,p_1,q_1),x_1-x_3\ra+\lm^2\|f(x_1,p_1,q_1)-f(x_3,p_1,q_1)\|^2\right]^\frac{1}{2}\\
&\qquad\qquad\qquad+\lm\|f(x_3,p_2,q_2)-f(x_3,p_1,q_1)\|\\
&\disp\le\left[\|x_1-x_3\|^2-2\lm\sigma\|x_1-x_3\|^2+\lm^2L^2\|x_1-x_3\|^2\right]^\frac{1}{2}+\lm L\big(d_1(p_1,p_2)+d_2(q_1,q_2)\big)\\
&\disp\le(1+\lm^2 L^2-2\lm\sigma)^\frac{1}{2}\|x_1-x_3\|+\lm L\sqrt{2\eta}d_1(p_1,p_2)^{\frac{1}{2}}+\lm Ld_2(q_1,q_2),\\
\end{array}\]
where the second inequality is due to \eqref{lipf}. Combining the last inequality with \eqref{e.7.0a} gives us that
\[
\|x_1-x_3\|\le\gg_1 d_1(p_1,p_2)^\frac{1}{2}+\gg_2 d_2(q_1,q_2)
\]
with the positive constants $\gg_1$ and $\gg_2$ defined by
\begin{eqnarray}\label{gg}
\gg_1:=\left(1-\frac{\sqrt{1+\lm^2L^2-2\lm\kk}}{1-r\lm}\right)^{-1}\left(\frac{\lm L\sqrt{2\eta}}{1-r\lm}+\ell\right),\gg_2:=\left(1-\frac{\sqrt{1+\lm^2L^2-2\lm\kk}}{1-r\lm}\right)^{-1}\frac{\lm L}{1-r\lm}.
\end{eqnarray}
Unifying it with estimate \eqref{e.7} ensures that
\[\begin{array}{ll}
\|v_1-v_2-2\kk(x_1-x_2)\|&\disp\le\|v_1-v_2-2\kk(x_3-x_2)\|+2\kk\|x_3-x_1\|\\
&\disp\le\|v_1-v_2\|+2\kk\gg_1 d_1(p_1,p_2)^\frac{1}{2}+2\kk\gg_2 d_2(q_1,q_2),
\end{array}\]
which verifies \eqref{4.7} and thus completes the proof of the theorem. \endproof\vspace*{-0.15in}

\begin{Example}{\bf (condition $\sigma>{\cal R}$ is essential for full stability).}\label{rm3}
{\rm Consider PVS \eqref{VS} with
\begin{eqnarray}\label{s-ex}
g(x,p):=\delta_{\R^2_+}(x)-\frac{r}{2}\|x\|^2\;\mbox{and}\;f(x,p,q):=\sigma x+p+q,\;(x,v,p,q)\in\R^2\times\R^2\times\R^2\times\R^2,
\end{eqnarray}
where $\sigma\le\mathcal{R}$. It is clear that $g$ is parametrically continuously prox-regular at $(\ox,\op)=(0,0)\in\R^2\times\R^2$ for $\ov=0\in\R^2$ with a prox-parameter $r>\mathcal{R}$ in \eqref{par-prox} and that assumptions {\rm(A1)--(A3)} are satisfied. Moreover, inequality \eqref{smon} holds and $x\in S(v,p,q)$ if and only if $v+(r-\sigma)x-p-q\in N_{\R^2_+}(x)$ that amounts to $(\sigma-r)x+p+q-v\in\R^2_+$ and $x\in\R^2_+$. Since $\sigma-r\le 0$, we get from the latter that $p+q-v\in\R^2_+$, which is not the case for all $(v,p,q)$ around $(0,0,0)\in\R^2\times\R^2\times\R^2$. Thus the solution map \eqref{ss} with the initial data \eqref{s-ex} does not admit a single-valued localization satisfying the H\"olderian property \eqref{4.7} around $(0,0,0)\in\R^2\times\R^2\times\R^2$ for $\ox$.}
\end{Example}\vspace*{-0.25in}

\begin{Remark}{\bf (on strong monotonicity).}\label{rem-mon} {\rm The strong monotonicity condition \eqref{smon}, together with its global and set-valued counterparts, have been well-recognized in variational analysis and its applications. It was used for sensitivity analysis of variational inequalities in the pioneering work by Dafermos \cite{DS} and then was profoundly developed by Yen \cite{Y2,Y1} in the study of Lipschitzian and H\"olderian properties of solutions to parametric variational and quasi-variational inequalities. Quite recently, Mordukhovich and Nghia \cite{MN2} established complete {\em coderivative characterizations} of this property and its set-valued counterpart in both finite and infinite dimensions.}
\end{Remark}\vspace*{-0.25in}

\begin{Remark}{\bf (H\"olderian properties of solutions to quasi-variational inequalities).}\label{ex1} {\rm Let $C\colon\mathcal{P}\tto X$ be a set-valued mapping with closed and convex values $C(p)$ for any $p\in\mathcal{P}$, and let $(\ox,\op)\in\gph C$ with $\hat v=\ov-f(\ox,\op,\oq)\in N_{C(\op)}(\ox)$. It is easy to see that the indicator function $g(x,p)=\delta_{C(p)}(x)$ is parametrically continuously prox-regular at any $(\ox,\op)\in\gph C$ for $\hat v$ with the threshold ${\cal R}=0$. The variational system \eqref{VS} reads now as
\begin{eqnarray}\label{e.7.7.1}
v\in f(x,p,q)+N_{C(p)}(x)
\end{eqnarray}
and belongs to the area of {\em quasi-variational inequalities}, which distinct from classical variational inequalities by the presence of parameters in convex sets $C(p)$. Considering \eqref{e.7.7.1} with $f=f(x,q)$ and imposing the local strong monotonicity of $f$ in $x$ \eqref{smon} together with the Lipschitz-like property of the mapping $C\colon\mathcal{P}\tto X$ (an equivalent of our (A3) assumption in this setting), Yen \cite[Theorem~2.1]{Y1} established the classical H\"older continuity of a single-valued localization of the solution map to \eqref{e.7.7.1}. Recall that the H\"olderian full stability property for \eqref{e.7.7.1}---and for the essentially more general setting of PVS \eqref{VS} treated in Theorem~\ref{thm2}---is {\em significantly stronger} than the standard H\"older continuity of localizations; see  more discussions and examples in \cite{MN2}.}
\end{Remark}\vspace*{-0.1in}

The next result is a Lipschitzian counterpart of Theorem~\eqref{thm2} for the nonsmooth PVS \eqref{VS}.\vspace*{-0.1in}

\begin{Theorem}{\bf (Lipschitzian full stability under strong monotonicity).}\label{thm2.1} Consider the setting of Theorem~{\rm\ref{thm2}} and assume in addition that the Lipschitz-like property in \eqref{K} is satisfied. Then $\ox$ is a Lipschitzian fully stable solution to PVS \eqref{VS} corresponding to the parameter triple $(\ov,\op,\oq)$.
\end{Theorem}\vspace*{-0.05in}
{\bf Proof.} It follows the lines in the proof of Theorem~\ref{thm2} with some modifications. Applying Theorem~\ref{prop2}, we find a single-valued localization $\pi_\lm$ of the proximal mapping $P_\lm$ at $(\lm\hat v+\ox,\op)$ for $\ox$, which satisfies the Lipschitzian estimate \eqref{e.7.1}. Thus there is a neighborhood $U_1\times P_1\times V_1\subset U\times\mathcal{P}\times X$ of $(\ox,\op,\tilde v)$ with $\tilde v:=\lm\hat v+\ox$ and some constants $\kk_0,\ell_0>0$  such that $\gph\pi_\lm=\gph P_\lm\cap(V_1\times P_1\times U_1)$ and \eqref{e.7.1} holds for all $(v_1,p_1),(v_2,p_2)\in V_1\times P_1$. This in turn yields
\begin{eqnarray}\label{e.7.4.1}
\|\pi_\lm(v_1,p_1)-\pi_\lm(v_2,p_2)\|\le\frac{1}{\kk_0}\|v_1-v_2\|+\frac{\ell_0}{2\kk_0}d_1(p_1,p_2).
\end{eqnarray}
Note furthermore that $x\in S(v,p,q)$ if and only if
\[
x+\lm v-\lm f(x,p,q)\in\lm\partial_x g(x,p)+x,
\]
which amounts to saying that $x\in P_\lm(x+\lm v-\lm f(x,p,q),p)$. By the Lipschitz continuity of $f$ in (A1), we find a neighborhood $U_2\times V_2\times P_2\times Q_2\subset U_{1}\times V_1\times P_{1}\times Q_{1}$ of $(\ox,\ov,\op,\oq)$ such that $(x+\lm v-\lm f(x,p,q),p)\in V_1\times P_1$ for all $(x,v,p,q)\in U_2\times V_2\times P_2\times Q_2$. Thus $x\in S(v,p,q)$ if and only if $x=\pi_\lm(x+\lm v-\lm f(x,p,q),p):=H_{(v,p,q)}(x)$ for $(x,v,p,q)\in U_2\times V_2\times P_2\times Q_2$. Fixing now $(v,p,q)\in V_2\times P_2\times Q_2$ and arguing similarly to the proof of Theorem~\ref{thm2}, we deduce that the mapping $H_{(v,p,q)}$ from \eqref{H} satisfies the contraction condition for some $\lm\in(0,\mathcal{R}^{-1})$, and therefore $H_{(v,p,q)}$ has a unique fixed point in $\B_\delta(\ox)$ denoted by $\vt(v,p,q)$.

It remains to verify the error bound estimate \eqref{mp1} of Lipschitzian full stability. Take any $(v_1,p_1),(v_2,p_2)\in\B_\eta(\ov,\op)$ with $x_1:=\vt(v_1,p_1,q_1)$, $x_2:=\vt(v_2,p_2,q_2)$, and $x_3:=\vt(v_1,p_2,q_2)$ and observe that $\gph\vt=\gph S\cap(\B_\eta(\ov)\times\B_\eta(\op)\times \B_\eta(\oq)\times\B_\delta(\ox))$. It gives us the representations
$$
x_i=\pi_\lm(x_i+\lm v_i-\lm f(x_i,p_i,q_i),p_i),\;i=1,2,\;\mbox{ and }\;x_3=\pi_\lm(x_3+\lm v_1-\lm f(x_1,p_2,q_2),p_2).
$$
Similar to the proof of Theorem~\ref{thm2}, we obtain the estimate
\begin{eqnarray}\label{e.7.6}
\|v_1-v_2-2\kk(x_3-x_2)\|\le\|v_1-v_2\|\;\mbox{with}\;\kk:=\lm^{-1}(\kk_0+\lm\sigma-1)=\lm^{-1}(\lm\sigma-\lm r)=\sigma-r.
\end{eqnarray}
It follows further from the Lipschitz continuity of $\pi_\lm$ in \eqref{e.7.4.1} that
\begin{equation}\label{e.7.7}
\begin{array}{ll}
\|x_1-x_3\|&\disp=\|\pi_\lm(x_1+\lm v_1-\lm f(x_1,p_1,q_1),p_1)-\pi_\lm(x_3+\lm v_1-\lm f(x_3,p_2,q_2),p_2)\|\\
&\disp\le\frac{1}{1-r\lm}\|x_1-x_3-\lm(f(x_1,p_1,q_1)-f(x_3,p_2,q_2))\|+\ell d_1(p_1,p_2)
\end{array}
\end{equation}
with $\ell=\frac{\ell_0}{2\kk_0}$, and therefore we get the inequalities
\[\begin{array}{ll}
&\|x_1-x_3-\lm(f(x_1,p_1,q_1)-f(x_3,p_2,q_2))\|\\
&\disp\le\|x_1-x_3-\lm (f(x_1,p_1,q_1)-f(x_3,p_1,q_1))\|+\lm \|f(x_3,p_2,q_2)-f(x_3,p_1,q_1)\|\\
&\disp\le\left[\|x_1-x_3\|^2-2\lm\la f(x_1,p_1,q_1)-f(x_3,p_1,q_1),x_1-x_3\ra+\lm^2\|f(x_1,p_1,q_1)-f(x_3,p_1,q_1)\|^2\right]^\frac{1}{2}\\
&\disp\qquad\qquad+\lm L\big(d_1(p_1,p_2)+d_2(q_1,q_2)\big)\\
&\disp\le(1+\lm^2 L^2-2\lm\sigma)^\frac{1}{2}\|x_1-x_3\|+\lm L\big(d_1(p_1,p_2)+d_2(q_1,q_2)\big),\\
\end{array}\]
where the second one comes from \eqref{lipf}. Combining the last inequality with \eqref{e.7.7} yields
\[
\|x_1-x_3\|\le\gg_1 d_1(p_1,p_2)+\gg_2 d_2(q_1,q_2),
\]
where the constant $\gg_2$ is defined in \eqref{gg} while
$$
\gg_1:=\left(1-\frac{\sqrt{1+\lm^2L^2-2\lm\sigma}}{1-r\lm}\right)^{-1}\left(\frac{\lm L}{1-r\lm}+\ell\right).
$$
This together with \eqref{e.7.6} ensures the estimates
\[\begin{array}{ll}
\|v_1-v_2-2\kk(x_1-x_2)\|&\disp\le \|v_1-v_2-2\kk(x_3-x_2)\|+2\kk\|x_3-x_1\|\\
&\disp \le \|v_1-v_2\|+2\kk\gg_1d_1(p_1,p_2)+2\kk\gg_2d_2(q_1,q_2),
\end{array}\]
which clearly verifies \eqref{mp1} and thus completes the proof of the theorem.\endproof

Finally in this section, we present an efficient consequence of Theorem~\ref{thm2.1} for the case of parametric {\em affine quasi-variational inequalities} formulated as:
\begin{eqnarray}\label{aqvi}
\begin{array}{ll}
&\mbox{find }\;x\in C(p):=\big\{x\in\R^n\big|\;Ax\le p,\;x\ge 0\big\},\quad p\in\R^l,\\
&\mbox{such that }\;\la f(x,p,q)-v,y-x\ra\le 0\;\mbox{ for all }\;y\in C(p),
\end{array}
\end{eqnarray}
where $A$ is an $d\times n$ matrix, $v\in\R^n$, $q\in\R^m$, and $f\colon\R^n\times\R^l\times\R^m\to\R^n$. It is clear that \eqref{aqvi} can be equivalently written in the generalized equation form \eqref{e.7.7.1} with $C(p)$ from \eqref{aqvi}.\vspace*{-0.06in}

\begin{Corollary}{\bf (Lipschitzian full stability for affine quasi-variational inequalities).}\label{lip-qvi} Let $g(x,p)=\delta_{C(p)}(x)$ in \eqref{VS}, where the sets $C(p)$ are defined in \eqref{aqvi}, and where the base mapping $f$ satisfies the standing assumption {\rm(A1)} and the local strong monotonicity condition \eqref{smon} around $(\ox,\op,\oq)$ with some modulus $\sigma>0$. Then $\ox\in S(\ov,\op,\oq)$ in \eqref{ss} is a Lipschitzian fully stable solution to \eqref{aqvi} corresponding to the parameter triple $(\ov,\op,\oq)$.
\end{Corollary}\vspace*{-0.05in}
{\bf Proof.} As mentioned in a more general setting of Remark~\ref{ex1}, the indicator function $g(x,p)=\delta_{C(p)}(x)$ is parametrically continuously prox-regular at any $(\ox,\op)\in\gph C$ for $\hat v=\ov-f(\ox,\op,\oq)$ with the threshold ${\cal R}=0$. Furthermore, the imposed assumptions (A3) and \eqref{K} of Theorem~\ref{thm2.1} reduce to the Lipschitz-like property of the mappings $p\mapsto\epi\dd_{C(p)}(\cdot)$ and $p\mapsto\gph N_{C(p)}(\cdot)$ around the point $(\op,(\ox, \hat v))$, respectively. Both these conditions hold automatically due to the polyhedral structure of $C(p)$; it follows from the classical results by Walkup and Wets \cite{ww}. Thus we deduce the claimed Lipschitzian full stability in \eqref{aqvi} from Theorem~\ref{thm2.1}. \endproof

Corollary~\ref{lip-qvi} extends the result by Yen \cite[Theorem~3.1]{Y2} for \eqref{aqvi} with $f=f(x,q)$ who established the existence of a Lipschitz continuous localization of the solution map to \eqref{aqvi}, which is a weaker property than the Lipschitzian full stability in Corollary~\ref{lip-qvi}. Note also that system \eqref{aqvi} belongs to the class of parametric {\em variational conditions} studied in Section~9 in {\em nonlinear} settings where certain {\em constraint qualifications} are needed to ensure the validity of (A3) and \eqref{K}.\vspace*{-0.2in}

\section{Full Stability in Parametric Variational Inequalities}
\setcounter{equation}{0}

In this section we consider a special subclass of PVS \eqref{VS} given as
\begin{eqnarray}\label{Ge}
v\in f(x,p)+N_C(x),
\end{eqnarray}
where $X$ is a Hilbert space, $(\mathcal{P},d)$ is a metric space, and where $C$ is a closed and {\em convex} subset of $X$. This model has been well recognized as Robinson's generalized equation form \cite{R0} of {\em parametric variational inequalities} (PVI). Indeed, similarly to \eqref{aqvi} it can be written in the classical variational inequality form, where the set $C$ is now parameter independent. The parametric sets of solutions to the parameterized variational inequality \eqref{Ge} are denoted by
\begin{eqnarray}\label{sol}
\Hat S(v,p):=\big\{x\in X\big|\;v\in f(x,p)+N_C(x)\big\}.
\end{eqnarray}

We concentrate in what follows on the study of {\em Lipschitzian} full stability of $\ox\in\Hat S(\ov,\op)$ in \eqref{Ge} while using the term ``full stability" for it. Observing that assumptions (A2) and (A3) on the potential $g(x)=\dd_C(x)$ are trivially satisfied, we keep here assumptions (A1) and (A4) on the base $f=f(x,p)$, which are the {\em standing assumptions} in this section.

Under these assumptions Section~3, which follows \cite{MN2}, contains second-order characterizations of full stability for general PVS \eqref{VS} presented in two forms: neighborhood in the case of Hilbert decision spaces and pointbased in finite dimensions. The main goal of this section is to derive {\em pointbased} conditions for full stability in PVI \eqref{Ge} in {\em infinite dimensions} provided that the set $C$ therein in {\em polyhedric}. The latter assumption has been well recognized and employed in the study and applications of elliptic partial differential equations.\vspace*{0.03in}

We start with two basic definitions for the major results of this section.\vspace*{-0.1in}

\begin{Definition}{\bf (Legendre forms).}\label{legen} The real-valued function $Q:X\to\R$ is a {\sc Legendre form} if it is weakly lower semicontinuous, represented as $Q(x)=\la Ax,x\ra$ with some linear operator $A:X\to X$, and satisfies the implication
$$
\big[x_k\st{w}\to x,\;Q(x_k)\to Q(x)\big]\Longrightarrow x_k\to x\;\mbox{ as }\;k\to\infty.
$$
\end{Definition}\vspace*{-0.05in}

To the best of our knowledge, this notion, which is automatic in finite dimensions, was introduced by Ioffe and Tikhomirov \cite[Section~6.2]{IT} and then was specified, developed, and  largely applied in, e.g., \cite{BBS,B,BS,MN} among many other publications. The next concept of polyhedricity was independently initiated by Haraux \cite{H} and Mignot \cite{Mi} motivated by applications to semilinear elliptic PDEs. The class of polyhedric sets includes all the convex polyhedra while being significantly broader than the latter, especially in infinite dimensions; see, e.g., \cite{BBS,B,BS,H,hms,HS,Mi,MN} for important examples, further discussions, and various applications.\vspace*{-0.1in}

\begin{Definition} {\bf (polyhedric sets).}\label{poly} Let $C$ be a closed and convex subset of $X$. It is said to be {\sc polyhedric} at $\ox\in C$ for some $\hat v\in N_C(\ox)$ if we have the representation
\begin{eqnarray}\label{5.3}
\K_C(\ox,\hat v):=T_C(\ox)\cap\{\hat v\}^\perp={\rm cl}\Big\{\mathcal{R}_C(\ox)\cap\{\hat v\}^\perp\Big\}
\end{eqnarray}
of the corresponding critical cone $\K(\ox,\hat v)$, where
\begin{equation}\label{5.3a}
\mathcal{R}_C(\ox):=\bigcup_{t>0}\Big[\frac{C-\ox}{t}\Big]
\end{equation}
is called the radial cone, and where $T_C(\ox):={\rm cl}\,\mathcal{R}_C(\ox)$ is the classical tangent cone to $C$ at $\ox$. If $C$ is polyhedric at each $\ox\in C$ for any $\hat v\in N_C(\ox)$, we say that the set $C$ is polyhedric.
\end{Definition}\vspace*{0.03in}

The following proposition, taken from \cite[Theorem~6.2]{MN} and employed below, provides a precise calculation of the regular coderivative for the normal cone mapping generated by a polyhedric set.\vspace*{-0.1in}

\begin{Proposition}{\bf (regular coderivative calculation for the normal cone mapping to polyhedric sets).}\label{pro2.3} For any $\ox\in C$ and $\hat v\in N_C(\ox)$ we have the inclusion
\begin{eqnarray*}
\dom\hat D^*N_C(\ox,\hat v)\subset-\K_C(\ox,\hat v)
\end{eqnarray*}
via the critical cone \eqref{5.3}. If in addition $C$ is polyhedric at $\ox\in C$ for $\hat v$, then
\begin{eqnarray}\label{5.5}
\Hat D^*N_C(\ox,\hat v)(w)=\K_C(\ox,\hat v)^*\quad\mbox{whenever}\quad w\in-\K_C(\ox,\hat v).
\end{eqnarray}
\end{Proposition}

Now we are ready to derive a major result of this section that gives a pointbased characterization of full stability for solutions to PVI \eqref{Ge} in Hilbert spaces. It can be seen as a far-going extension of \cite[Theorem~7.2]{MN}, which characterizes Lipschitzian full stability of local minimizers in optimal control problems governed by semilinear PDEs with Legendre elliptic operators.\vspace*{-0.1in}

\begin{Theorem}{\bf (pointbased characterization of full stability for PVI solutions).}\label{NGE} Let $\ox\in\Hat S(\ov,\op)$ in \eqref{sol} with $\Hat v:=\ov-f(\ox,\op)$, and let the base function $f$ satisfy the standing assumptions of this section. Consider the following statements:
	
{\bf (i)} $\ox$ is a fully stable solution to PVI \eqref{Ge} corresponding to the parameter pair $(\ov,\op)$.
	
{\bf (ii)} We have the positive-definiteness condition
\begin{eqnarray}\label{5.7}
\la\nabla_x f(\ox,\op)w,w\ra>0\quad\mbox{for all}\quad w\in{\cal H}^{\rm w}(\ox,\hat v),\;w\ne 0
\end{eqnarray}
in terms of the $($weak$)$ sequential outer limit \eqref{pk} of the critical cones
\begin{eqnarray}\label{5.8}
{\cal H}^{\rm w}(\ox,\hat v):=\Limsup_{(x,v)\st{{\rm gph}\,N_C}\longrightarrow(\ox,\hat v)}\K_C(x,v).
\end{eqnarray}
	
{\bf (iii)} We have another positive-definiteness condition
\begin{eqnarray}\label{5.9}
\la\nabla_x f(\ox,\op)w,w\ra>0\quad\mbox{for all}\quad w\in{\cal H}^{\rm s}(\ox,\hat v),\;w\ne 0,
\end{eqnarray}
where ${\cal H}^{\rm s}(\ox,\hat v)$ stands for the strong counterpart of \eqref{pk} via the norm topology on $X$ defined by
\begin{eqnarray}\label{5.10}
\begin{array}{ll}
{\cal H}^{\rm s}(\ox,\hat v)&\disp:={\rm s}-\Limsup_{(x,v)\st{{\rm gph}\,N_C}\longrightarrow(\ox,\hat v)}\K_C(x,v)\\
&\disp=\Big\{z\in X\big|\;\exists\;{\rm sequences }\,(x_k,v_k)\st{{\rm gph}N_C}\longrightarrow(\ox,\hat v),\;z_k\in \K_C(x_k,v_k),\;\|z_k-z\|\to 0\Big\}.
\end{array}
\end{eqnarray}
Then the following assertions are satisfied:
	
{\bf(A)} If $Q(w):=\la\nabla_x f(\ox,\op)w,w\ra$ is Legendre, then {\rm(ii)}$\Longrightarrow${\rm (i)}.
	
{\bf(B)} If the $C$ is polyhedric, then {\rm(i)}$\Longrightarrow${\rm(iii)}.
\end{Theorem}\vspace*{-0.07in}
{\bf Proof.} To justify assertion (A), we apply Proposition~\ref{Lippq} to PVI \eqref{Ge} and observe that the only thing to verify is the validity of condition \eqref{4.8b}. Arguing by contradiction, suppose that the latter condition fails and then find a sequence $(u_k,v_k,w_k,z_k)\in X\times X\times X\times X$ such that $(u_k,v_k)\st{{\rm gph}\,N_C}\to(\ox,\hat v)$, $z_k\in\Hat D^*N_C(u_k,v_k)(w_k)$, and
\begin{eqnarray}\label{5.11}
\la\nabla_x f(\ox,\op)w_k,w_k\ra+\la z_k,w_k\ra<\frac{1}{k}\|w_k\|^2\;\mbox{ for all }\;k\in\IN.
\end{eqnarray}
Since the normal cone operator $N_C\colon X\tto X$ is maximal monotone by the classical Rockafellar's theorem due to the convexity of $C$, it follows from \cite[Lemma~3.3]{CT} that $\la z_k,w_k\ra\ge 0$. Further, we get from Proposition~\ref{pro2.3} that $w_k\in-{\cal K}_C(u_k,v_k)$. Combining these with \eqref{5.11} ensures that
\begin{eqnarray}\label{Q0}
\la\nabla_x f(\ox,\op)w_k,w_k\ra<\frac{1}{k}\|w_k\|^2\quad\mbox{with}\quad w_k\in-{\cal K}_C(u_k,v_k),
\end{eqnarray}
which yields $w_k\ne 0$. Defining $\ow_k:=w_k\|w_k\|^{-1}$ with $\|\ow_k\|=1$, we deduce from \eqref{Q0} that
\begin{eqnarray}\label{Q}
Q(\ow_k)=\la\nabla_x f(\ox,\op)\ow_k,\ow_k\ra<\frac{1}{k}\quad\mbox{with}\quad\ow_k\in-{\cal K}_C(u_k,v_k),
\end{eqnarray}
and so there exists a subsequence of $\{\ow_k\}$ (without relabeling) that weakly converges to some $\ow$. By $(u_k,v_k)\st{{\rm gph}\,N_C}\longrightarrow(\ox,\hat v)$ it follows from \eqref{5.8} that $\ow\in-{\cal H}^{\rm w}(\ox,\hat v)$. Employing the imposed weak lower semicontinuity of $Q$ tells us by \eqref{Q} and \eqref{5.7} that
\begin{eqnarray}\label{Q1}
0\le Q(\ow)\le\liminf_{k\to\infty}Q(\ow_k)\le 0,
\end{eqnarray}
and therefore $Q(\ow)=0$. Since $Q$ is Legendre, it follows from \eqref{Q1} and Definition~\ref{legen} that $\|\ow_k-\ow\|\to 0$ along a subsequence of $k\to\infty$, which gives us $\|\ow\|=1$, $\ow\in{\cal H}^{\rm w}(\ox,\hat v)$, and $Q(\ow)=Q(-\ow)=0$. This contradicts \eqref{5.7} and hence completes the proof of assertion (A).\vspace*{0.03in}

Now we verify assertion (B) assuming that $C$ is polyhedric. It follows from full stability in (i) and Proposition~\ref{Lippq} in the case of $g(x,q)=\delta_C(x)$ that there are numbers $\kk,\eta>0$ such that for any $(u,v)\in\gph N_C\cap\B_\eta(\ox,\hat v)$ we have the condition
\begin{eqnarray}\label{5.13}
\la\nabla_x f(\ox,\op)w,w\ra+\la z,w\ra\ge\kk\|w\|^2\quad\mbox{whenever}\quad z\in\Hat D^*N_C(u,v)(w),w\in X.
\end{eqnarray}
Since $C$ is polyhedric, Proposition~\ref{pro2.3} tells us $w\in-{\cal K}(u,v)$ and $\Hat D^*N_C(u,v)(w)={\cal K}_C(u,v)^*$, which yield $0\in\Hat D^*N_C(u,v)(w)$. This together with \eqref{5.13} and \eqref{5.5} shows that
\begin{eqnarray}\label{5.13a}
\la\nabla_x f(\ox,\op)w,w\ra\ge\kk\|w\|^2\quad\mbox{for all}\quad w\in-{\cal K}_C(u,v)
\end{eqnarray}
when $(u,v)\in\gph N_C\cap\B_\eta(\ox,\hat v)$. Passing to the limit in \eqref{5.13a} as $\eta\dn 0$ and using the strong convergence in \eqref{5.10}, we arrive at \eqref{5.9} and hence finish the proof of the theorem.\endproof

The next lemma effectively estimates the limiting forms ${\cal H}^{\rm w}$ and ${\cal H}^{\rm s}$ in Theorem~\ref{NGE} via the tangent and critical cones for the set $C$, which allows us to establish more direct and verifiable conditions for full stability of solutions to \eqref{Ge} in both finite and infinite dimensions.\vspace*{-0.1in}

\begin{Lemma}{\bf(estimates of weak and strong outer limits of the critical cone).}\label{pro5.5} In the general setting of Theorem~{\rm\ref{NGE}} we have the upper estimate
\begin{eqnarray}\label{sub}
{\cal H}^{\rm w}(\ox,\hat v)\subset{\rm cl}\big[T_C(\ox)-T_C(\ox)\big]\cap\{\hat v\}^\perp.
\end{eqnarray}
If in addition $C$ is polyhedric, then the lower estimate
\begin{eqnarray}\label{sup}
{\cal K}_C(\ox,\hat v)-{\cal K}_C(\ox,\hat v)\subset{\cal H}^{\rm s}(\ox,\hat v)
\end{eqnarray}
is satisfied. It furthermore $\dim X<\infty$ and the set $C$ is polyhedral, then we have the representation
\begin{eqnarray}\label{5.16}
{\cal H}^{\rm s}(\ox,\hat v)={\cal H}^{\rm w}(\ox,\hat v)={\cal K}_C(\ox,\hat v)-{\cal K}_C(\ox,\hat v).
\end{eqnarray}
\end{Lemma}\vspace*{-0.05in}
{\bf Proof.} To verify \eqref{sub}, pick $w\in{\cal H}^{\rm w}(\ox,\hat v)$ and find a sequence $(u_k,v_k,w_k)$ such that $w_k\st{w}\to w$, $(u_k,v_k)\st{{\rm gph}N_C}\longrightarrow\ox,\hat v)$, and $w_k\in{\cal K}_C(u_k,v_k)=T_C(u_k)\cap\{v_k\}^\perp$. Hence for each $k\in\IN$ there are sequences $C\ni y_{n_k}\to u_k$ and $t_{n_k}\dn 0$ with $\frac{y_{n_k}-u_k}{t_{n_k}}\to w_k$ as $n\to\infty$. In this way we construct
$$
(z_k,\al_k)\in\big\{(y_{n_k},t_{n_k})\big|\;n\in\IN\big\}
$$
so that $z_k\to\ox$, $\al_k\dn 0$ , and $\frac{z_{k}-u_k}{\al_k}-w_k\to 0$ as $k\to\infty$. It follows that
\[
w={\rm w}-\lim_{k\to\infty}w_k={\rm w}-\lim_{k\to\infty}\frac{z_{k}-u_k}{\al_k}={\rm w}-\lim_{k\to\infty}\frac{(z_{k}-\ox)-(u_k-\ox)}{\al_k}\in {\rm cl}^{\rm w}\,\Big[T_C(\ox)-T_C(\ox)\Big],
\]
where the symbol ``${\rm w}-\lim$" indicates that the weak limit in $X$ is taken. Since $\la w_k,v_k\ra=0$, we get $\la w,\hat v\ra=0$ by passing to the limit, and hence
$$
w\in {\rm cl}^{\rm w}\big[T_C(\ox)-T_C(\ox)\big]\cap\{\hat v\}^\perp.
$$
Note that $T_C(\ox)-T_C(\ox)$ is a convex set. The classical Mazur theorem tells us that ${\rm cl}^{\rm w}\big[T_C(\ox)-T_C(\ox)\big]={\rm cl}\big[T_C(\ox)-T_C(\ox)\big]$. Therefore $w\in{\rm cl}\big[T_C(\ox)-T_C(\ox)\big]\cap\{\hat v\}^\perp$, which justifies \eqref{sub}.\vspace*{0.03in}

Now suppose that $C$ is polyhedric. To verify \eqref{sup}, pick any $w=w_1-w_2$ from the left-hand side of \eqref{sup} with $w_1,w_2\in {\cal K}_C(\ox,\hat v)$. The polyhedricity of $C$ allows us to find \eqref{5.3} sequences $w_{1k}\to w_1$, $w_{2k}\to w_2$, and $t_{1k},t_{2k}\dn 0$ such that $\ox+t_{1k}w_{1k}\in C$, $\ox+t_{2k}w_{2k}\in C$, and $w_{1k}, w_{2k}\in\{\hat v\}^\perp$. Defining $t_k:=\min\{t_{1k},t_{2k}\}$, we get from the convexity of $C$ that $u_k:=\ox+t_kw_{1k}=(1-t_kt_{1k}^{-1})\ox+t_kt_{1k}^{-1}(\ox+t_{1k}w_{1k})\in C$ and similarly $x_k:=\ox+t_kw_{2k}\in C$. Thus it follows from \eqref{5.3} and \eqref{5.3a} that
\[
w_{2k}-w_{1k}=\frac{x_k-u_k}{t_k}\in{\cal R}_C(u_k)\cap\{\hat v\}^\perp\subset\K_C(u_k,\hat v).
\]
Since $w_{2k}-w_{1k}\to w_2-w_1=w$ and $u_k\to\ox$ as $k\to\infty$, we deduce from \eqref{5.10} that $w\in{\cal H}^{\rm s}(\ox,\hat v)$ and hence complete the proof of inclusion \eqref{sup}.

Finally, assuming $\dim X<\infty$ ensures that ${\cal H}^{\rm w}(\ox,\hat v)={\cal H}^{\rm s}(\ox,\hat v)$. To verify \eqref{5.16}, it remains to show  by \eqref{sup} that the opposite inclusion holds therein when $C$ is polyhedral. Picking any $w\in{\cal H}^{\rm s}(\ox,\hat v)$, find a sequence $(w_k, x_k,v_k)$ such that $w_k\in{\cal K}_C(x_k,v_k)$, $(x_k,v_k)\st{{\rm gph}N_C}\to(\ox,\hat v)$ and $w_k\to w$ as $k\to\infty$. It follows from the formulas for the tangent cone and for the normal cone to a polyhedral set in the proof of \cite[Theorem~2]{DR1} that
\begin{eqnarray}
&&w_k\in T_C(x_k)=T_C(\ox)+\R\{x_k-\ox\},\label{TC}\\
&&v_k\in N_C(x_k)=N_C(\ox)\cap\{x_k-\ox\}^\perp\label{NC}
\end{eqnarray}
for all large $k\in\IN$ and that $T_C(\ox)\cap\{v_k\}^\perp\subset T_C(\ox)\cap\{\hat v\}^\perp={\cal K}_C(\ox,\hat v)$. Since $\R_+\{x_k-\ox\}\subset T_C(\ox)$ and the cone $T_C(\ox)$ is convex, we obtain from \eqref{TC} that
\[
w_k\in T_C(x_k)=T_C(\ox)-\R_+\{x_k-\ox\}.
\]
Hence there exist $y_k\in T_C(\ox)$ and $r_k\ge 0$ with $w_k=y_k-r_k(x_k-\ox)$. Since $v_k\in N_C(x_k)$, we get $\la x_k-\ox,v_k\ra=0$ from \eqref{NC}. It follows furthermore that $0=\la w_k,v_k\ra=\la y_k,v_k\ra$, and so $y_k\in T_C(\ox)\cap\{v_k\}^\perp\subset{\cal K}_C(\ox,\hat v)$. Observing from \eqref{NC} that $r_k(x_k-\ox)\subset T_C(\ox)\cap\{v_k\}^\perp\subset{\cal K}(\ox,\hat v)$ and using the above representation of $w_k$ give us the inclusion
\[
w_k\in{\cal K}_C(\ox,\hat v)-{\cal K}_C(\ox,\hat v)
\]
showing that $w\in{\rm cl}\,\big[{\cal K}_C(\ox,\hat v)-{\cal K}_C(\ox,\hat v)\big]={\cal K}_C(\ox,\hat v)-{\cal K}_C(\ox,\hat v)$, where the equality holds since ${\cal K}_C(\ox,\hat v)-{\cal K}_C(\ox,\hat v)$ is a subspace in finite dimensions. This verifies that \eqref{5.16} is satisfied.\endproof

The obtained results lead us to the following verifiable conditions for full stability in \eqref{Ge}.\vspace*{-0.1in}

\begin{Theorem}{\bf(refined pointbased conditions for full stability of solutions to PVI).}\label{coro5.6} In the framework of Theorem~{\rm\ref{NGE}}, consider the statements:
	
{\bf (i)} $\ox\in\Hat S(\ov,\op)$ is a fully stable solution to \eqref{Ge}.
	
{\bf (ii)} We have the positive-definiteness condition of the closure type
\begin{eqnarray}\label{5.19}
\la\nabla_x f(\ox,\op)w,w\ra>0\quad\mbox{for all}\quad w\in{\rm cl}\big[T_C(\ox)-T_C(\ox)\big]\cap\{\hat v\}^\perp,\;w\ne 0.
\end{eqnarray}
	
{\bf (iii)} We have the the positive-definiteness condition via the critical cone
\begin{eqnarray}\label{5.20}
\la\nabla_x f(\ox,\op)w,w\ra>0\quad\mbox{for all}\quad w\in{\cal K}_C(\ox,\hat v)-{\cal K}_C(\ox,\hat v),\;w\ne 0.
\end{eqnarray}
Then the following assertions hold:
	
{\bf(A)} If $Q(w)=\la\nabla_x f(\ox,\op)w,w\ra$ is Legendre, then {\rm(ii)}$\Longrightarrow${\rm(i)}.
	
{\bf(B)} If the set $C$ is polyhedric, then {\rm(i)}$\Longrightarrow${\rm(iii)}.
	
{\bf(C)} If $\dim X<\infty$ and the set $C$ is polyhedral, then the conditions in {\rm(ii)} and {\rm(iii)} are equivalent while being necessary and sufficient for the validity of {\rm(i)}.
\end{Theorem}\vspace*{-0.07in}
{\bf Proof.} It follows by combining the corresponding assertions in Theorem~\ref{NGE} and Lemma~\ref{pro5.5}.\endproof

Note that the positive-definiteness condition of the closure type \eqref{5.19} has been used in \cite[p.\ 405]{BS} to ensure that the solution map $\Hat S$ in \eqref{sol} has a single-valued and Lipschitz continuous localization around $(\ov,\op)$ provided that $Q(w)$ is Legendre. As discussed above, our assertion {\rm(A)} of Theorem~\ref{coro5.6} establishes the much stronger property of (Lipschitzian) {\em full stability} of solutions to PVI \eqref{Ge} under the same assumptions as in \cite{BS}.\vspace*{0.03in}

Finally in this section, we focus on the space $X=L^2(\O)$, where $\O$ is an open subset of $\R^n$. Consider the closed and convex set $C\subset X$ in \eqref{Ge} defined by the {\em magnitude constraint system}
\begin{eqnarray}\label{C}
C=\big\{x\in L^2(\O)\big|\;a\le x(y)\le b\;\mbox{ a.e. }\;y\in\O\big\}
\end{eqnarray}
with $-\infty\le a<b\le\infty$. Such constraints are typical in applications to PDE control of elliptic equations; see, e.g., \cite{B,BS,HS,hms,MN}.
It follows from \cite[Proposition~6.33]{BS} that the set $C\subset L^2(\O)$ in \eqref{C} is polyhedric. The next proposition taken from \cite[Proposition~7.3]{MN} gives us precise calculations of the strong and weak outer limits of critical cones in \eqref{5.8} and \eqref{5.10}, respectively.\vspace*{-0.1in}

\begin{Proposition}{\bf(computation of outer limits of critical cones).}\label{lim-calc} Let $C\subset L^2(\O)$ be given in \eqref{C}, and let $(\ox,\hat v)\in\gph N_C$. Then the limiting sets in \eqref{5.8} and \eqref{5.10} are calculated by
\begin{eqnarray*}
{\cal H}^{\rm s}(\ox,\hat v)={\cal H}^{\rm w}(\ox,\hat v)=\big\{u\in L^2(\O)\big|\; u(y)\hat v(y)=0\;\;\mbox{a.e.}\;\; \mbox{on}\;\;\O\big\}.
\end{eqnarray*}
\end{Proposition}\vspace*{-0.05in}

This proposition allows us to derive from Theorem~\ref{NGE} the {\em pointbased characterization} of full stability for solutions to infinite-dimensional PVI \eqref{Ge} with constrained sets of type \eqref{C}.\vspace*{-0.1in}

\begin{Corollary}{\bf (characterization of full stability for PVI generated by magnitude constrained systems).} Let $X=L^2(\O)$ with $C$ given in  \eqref{C}, and let $\ox\in\hat S(\ov,\op)$ in \eqref{sol}. In addition to the standing assumptions of this section, suppose that $Q(w):=\la\nabla f(\ox,\op)w,w\ra$ for $w\in L^2(\O)$ is a Legendre form. Then the full stability $\ox$ in \eqref{Ge} is equivalent to the validity of the pointbased positive-definiteness condition
\begin{eqnarray*}
\la\nabla_x f(\ox,\op)w,w\ra>0\quad\mbox{for all}\quad w\in L^2(\O)\setminus\{0\}\quad\mbox{with}\quad w(y)\big(\ov(y)-f(\ox,\op)(y)\big)=0\;\; \mbox{a.e.}\;\;y\in\O.
\end{eqnarray*}
\end{Corollary}\vspace*{-0.05in}
{\bf Proof.} Follows directly from Theorem~\ref{NGE} by employing the calculations of Proposition~\ref{lim-calc}.\endproof\vspace*{-0.1in}

\section{Full Stability in Parametric Variational Conditions}
\setcounter{equation}{0}

Here we study a class of finite-dimensional PVS \eqref{VS} known as the {\em parametric variational conditions} (PVC) that are given in the form
\begin{eqnarray}\label{VC}
v\in f(x,p,q)+N_{C(p)}(x)\;\mbox{ with }\;x\in X=\R^n,\;p\in\mathcal{P}=\R^l,\;q\in\mathcal{Q}=\R^m,
\end{eqnarray}
and the sets $C(p)$ defined by the inequality constraints
\begin{eqnarray}\label{6.2}
C(p):=\big\{x\in\R^n\big|\;\ph_i(x,p)\le 0\;\mbox{ for }\;i=1,\ldots,s\big\}
\end{eqnarray}
via the functions $\ph_i:\R^n\times\R^l\to\R$, $i=1,\ldots,s$, which are ${\cal C}^2$-smooth around a feasible point $(\ox,\op)\in\R^n\times\R^l$ of \eqref{6.2}. As mentioned above, inclusion \eqref{VC} can be written in form \eqref{VS} with $g(x,p)=\delta_{C(p)}(x)$ for $(x,p)\in\R^n\times\R^l$. Furthermore, \eqref{VC} can be represented as a quasi-variational inequality when the sets $C(p)$ are convex; see \cite{mo} and Section~7 for more discussions. The parametric solution map to \eqref{VC} is denoted by
\begin{eqnarray}\label{6.3}
\breve S(v,p,q):=\big\{x\in\R^n\big|\;v\in f(x,p,q)+N_{C(p)}(x)\big\}.
\end{eqnarray}

The standing assumptions of the base $f$ in \eqref{VC} are the same (A1) and (A4) as in Section~8. Our main aim here is to derive {\em characterization} of (Lipschitzian) full stability for solutions to \eqref{VC} expressed entirely in terms of the {\em initial data} of \eqref{VC} and \eqref{6.2}.

Given $(\ox,\op)$ satisfying \eqref{6.2}, recall that the partial {\em Mangasarian-Fromovitz constraint qualification} (MFCQ) with respect to $x$ holds at $(\ox,\op)$ if there is $d\in X$ such that
\begin{eqnarray}\label{mfcq}
\la\nabla_x\ph_i(\ox,\op),d\ra<0\;\mbox{ for }\;i\in I(\ox,\op):=\big\{i\in\{1,\ldots,s\}\big|\;\ph_i(\ox,\op)=0\big\}.
\end{eqnarray}
The {\em Lagrangian function} for the variational system in \eqref{VC} and \eqref{6.2} is defined by
$$
L(x,p,q,\lm):=f(x,p,q)+\disp\sum_{i=1}^{s}\lm_i\nabla_x\ph_i(x,p)\;\mbox{ with }\;x\in\R^n,\;p\in\R^l,\;\mbox{ and }\;\lm\in\R^m.
$$
It is well known under the validity of the partial MFCQ condition \eqref{mfcq} that for all feasible $(x,p)$ around $(\ox,\op)$ we have the representation
\begin{eqnarray}\label{6.4}
\begin{array}{ll}
f(x,p,q)+N_{C(p)}(x)=\big\{L(x,p,q,\lm)\big|\;\lm\in N\big(\ph(x,p);\Th\big)\big\}:=\Psi(x,p,q)\;\mbox{ with }\;\Th:=\R^s_-
\end{array}
\end{eqnarray}
and $\ph=(\ph_1,\ldots,\ph_s)\colon\R^n\times\R^l\to\R^s$. Hence any vector $\ov\in f(\ox,\op,\oq)+N_{C(\op)}(\ox)$ satisfies
\begin{equation}\label{6.5}
\ov\in\nabla_xL(\ox,\op,\oq,\lm)\;\mbox{ with some }\;\lm\in N_\Theta\big(\ph(\ox,\op)\big),
\end{equation}
and the set of {\em Lagrange multipliers} is given by
\begin{eqnarray}\label{6.6}
\Lm(\ox,\op,\oq,\ov):=\big\{\lm\in\R_+^s\big|\;\ov\in\nabla_x L(\ox,\op,\oq,\lm),\;\la\lm,\ph(\ox,\op)\ra=0\big\}.
\end{eqnarray}
It follows from \cite[Proposition~2.2]{LPR} that MFCQ \eqref{mfcq} implies that BCQ \eqref{bcq1} holds for $g(x,p)=\delta_{C(p)}(x)$ and that $g$ is parametrically continuously prox-regular at $(\ox,\op)$ for $\hat v=\ov-f(\ox,\op,\oq)$.

Recall further the {\em general strong second-order sufficient condition} (GSSOSC) for the variational condition \eqref{VC} as formulated by Kyparisis \cite{K2}: given $(\ox,\op)$ satisfying \eqref{6.2} and given $\ov\in f(\ox,\op,\oq)+N_{C(\op)}(\ox)$, the GSSOSC holds at $(\ox,\op,\oq,\ov)$ if for all $\lm\in\Lm(\ox,\op,\oq,\ov)$ we have
\begin{eqnarray}\label{6.8}
\la u,\nabla_{x}L(\ox,\op,\oq,\lm)u\ra>0\;\mbox{ whenever }\;\la\nabla_x\ph_i(\ox,\op),u\ra=0\;\mbox{ as }\;i\in I_+(\ox,\op,\lm),\;u\ne 0
\end{eqnarray}
with the strict complementarity index set $I_+(\ox,\op,\lm):=\big\{i\in I(\ox,\op)\big|\;\lm_i>0\big\}$. This condition is a slight modification and adaptation to \eqref{VC}, \eqref{6.2} of the strong second-order sufficient condition introduced by Robinson \cite{Ro} for nonlinear programs with ${\mathcal C}^2$-smooth data; cf.\ also Kojima \cite{ko}.

Next we modify for the case of PVC in \eqref{VC} and \eqref{6.2} the uniform second-order sufficient condition introduced recently by Mordukhovich and Nghia \cite{MN} for parametric nonlinear programs.\vspace*{-0.1in}

\begin{Definition} {\bf (general uniform second-order sufficient condition).}\label{gusosc} We say that the {\sc general uniform second-order sufficient condition} {\rm (GUSOSC)} holds at $(\ox,\op,\oq)$ satisfying \eqref{6.2} with $\ov\in\Psi(\ox,\op,\oq)$ if there are positive numbers $\eta,\ell$ such that
\begin{eqnarray}\label{gus}
\begin{array}{ll}
\disp\la\nabla_{x}L(x,p,q,\lm)u,u\ra\ge\ell\|u\|^2\;\mbox{ for all }\;(x,p,q,v)\in\gph\Psi\cap\B_\eta(\ox,\op,\oq,\ov),\;\lm\in\Lm(x,p,q,v),&\\
\la\nabla_x\ph_i(x,p),u\ra=0\;\mbox{ as }\;i\in I_+(x,p,\lm)\;\mbox{ and }\;\la\nabla_x\ph_i(x,p),u\ra\ge 0\;\mbox{ as }\;i\in I(x,p)\setminus I_+(x,p,\lm),
\end{array}
\end{eqnarray}
where the mapping $\Psi$ and the set $\Lm(x,p,q,v)$ are defined in \eqref{6.4} and \eqref{6.6}, respectively.
\end{Definition}\vspace*{-0.05in}

Similarly to the proof of \cite[Proposition~4.2]{MN} in the nonparametric setting we can check that, under the validity of the partial MFCQ \eqref{mfcq}, the GSSOSC from \eqref{6.8} implies the GUSOSC from Definition~{\rm\ref{gusosc}} at $(\ox,\op,\oq)$ with $\ov\in\Psi(\ox,\op,\oq)$ by passing to the limit in \eqref{gus}.\vspace*{0.03in}

The last qualification condition needed below is the partial {\em constant rank constraint qualification} (CRCQ) formulated as follows; cf.\ \cite{FP,L2}. We say that the {\em partial CRCQ} with respect to $x$ holds at $(\ox,\op)$ feasible to \eqref{6.2} if there is a neighborhood $W$ of $(\ox,\op)$ such that for any $J\subset I(\ox,\op)$ the gradient family $\big\{\nabla_x\ph_i(x,p)\big|\;i\in J\big\}$ has the same rank in $W$. It occurs that the simultaneous fulfillment of the partial MFCQ and CRCQ ensures the Lipschitz-like property of the graphical mapping $K$ from \eqref{K} crucial for the Lipschitzian full stability results obtained in Sections~3, 6, and 7. Recall that the MFCQ condition is not need while CRCQ is automatics in the case of affine PVC/quasi-variational inequalities considered in Corollary~\ref{lip-qvi}.\vspace*{-0.1in}

\begin{Proposition}{\bf (graphical Lipschitz-like property under partial MFCQ and CRCQ).}\label{cr} Assume that both partial MFCQ and CRCQ conditions hold at the point $(\ox,\op)$ feasible to \eqref{6.2}. Then, given any vector $\ov$ satisfying \eqref{6.5}, the Lipschitz-like property \eqref{K} holds around $(\op,\ox,\hat v)$ with $g(x,p)=\delta_{C(p)}(x)$ and $\hat v=\ov-f(\ox,\op,\oq)$.
\end{Proposition}\vspace*{-0.07in}
{\bf Proof.} Follows directly from \cite[Proposition~5.2]{MN} for the case of constant cost functions.\endproof

The next theorem constitutes the main result of this section.\vspace*{-0.1in}

\begin{Theorem}{\bf (characterization of full stability for PVC under partial MFCQ and CRCQ).}\label{thm6.3} Take $(\ox,\ov,\op,\oq)\in\R^n\times\R^n\times\R^l\times \R^m$ with $\ov\in\breve S(\ox,\op)$ from \eqref{6.3} and suppose in addition to the standing assumptions that both partial MFCQ and CRCQ conditions formulated above hold at $(\ox,\op,\oq)$. Then the following assertions are equivalent:

{\bf (i)} $\ox$ is fully stable solution to PVC \eqref{VC}, \eqref{6.2} corresponding to $(\ov,\op,\oq)$.

{\bf(ii)} GUSOSC from Definition~{\rm\ref{gusosc}} holds at $(\ox,\op,\oq)$ with $\ov\in\Psi(\ox,\op,\oq)$.\\
\end{Theorem}\vspace*{-0.25in}
{\bf Proof.} It suffices to show by Proposition~\ref{Lippq} that the imposed GUSOSC is equivalent to the second-order subdifferential condition \eqref{4.43}. This can be done via calculating the term $\la z,w\ra$ in \eqref{4.43} by using the formula for $\Hat D^*\partial g_p$ established in \cite[Theorem~6]{HKO} under the validity of MFCQ and CRCQ. The proof is similar to the one given in \cite[Theorem~5.3]{MN}, and so we omit details. \endproof

Another sufficient condition for the existence of a single-valued and Lipschitz continuous localization of the solution map $\breve S$ around $(\ov,\op,\oq,\ox)$, named the {\em strong coherent orientation condition} (SCOC), was obtained by Facchinei and Pang \cite{FP} with imposing both MFCQ and CRCQ and assuming in addition that all the functions $\ph_i$ in \eqref{6.2} are convex. The latter convex assumption was dismissed in the more recent paper by Lu \cite{L2}. Observe that both developments in \cite{FP,L2} rely on topological degree theory the application of which to sensitivity analysis in optimization and related areas was initiated by Kojima \cite{ko}. The previous result in this direction was obtained by Kyparisis \cite{K2} who proved, based on the implicit function techniques by Robinson \cite{R}, the local single-valuedness and continuity (while not Lipschitz continuity) of the solution map $\breve S$ for convex PVC \eqref{VC} under the so-called ``general modified strong second-order condition" that is stronger than GSSOSC \eqref{6.8}, which in turn is stronger than SCOC.

The next example demonstrates that our new GUSOSC from Definition~\ref{gusosc} is {\em strictly weaker} than GSSOSC and is {\em not implied} by SCOC even in the case of constraint functions linear with respect to $x\in\R^3$ and $p\in\R^2$ as well as cost functions linear in $p$ and quadratic in $x$ under the validity of both MFCQ and CRCQ. This shows that Theorem~\ref{thm6.3}, which completely characterizes full stability of solutions to nonconvex PVC, provides in particular new sufficient conditions for the local single-valuedness and Lipschitz continuity of the solution map \eqref{6.3} that are independent of \cite{FP,L2} and significantly extend those obtained in \cite{K2}.\vspace*{-0.1in}

\begin{Example}{\bf (refined sufficient conditions for single-valuedness and Lipschitz continuity of PVC under partial MFCQ and CRCQ).}\label{ex}
{\rm Consider PVC in \eqref{VC} and \eqref{6.2} with
\begin{eqnarray}\label{5.30}
\left\{\begin{array}{ll}
f(x,p)=\nabla\ph_0(x,p)\;\mbox{ for }\;\ph_0(x,p):=x_3+\Big(\frac{1}{4}+p_2\Big)x_1+p_1x_2+x^2_3-x_1x_2,\\
\ph_1(x,p):=x_1-x_3-p_1\le 0,\\
\ph_2(x,p):=-x_1-x_3+p_1\le 0,\\
\ph_3(x,p):=x_2-x_3-p_2\le 0,\\
\ph_4(x,p):=-x_2-x_3+p_2\le 0,\\
x=(x_1,x_2,x_3)\in\R^3,\;p=(p_1,p_2)\in\R^2.
\end{array}\right.
\end{eqnarray}
We can easily see that both partial MFCQ and CRCQ hold at $(\ox,\op)$ with $\ox=(0,0,0)$, $\op=(0,0)$. It follows from the arguments in \cite[Example 6.4]{MN} that our GUSOSC holds in this setting. Thus we have by Theorem~\ref{thm6.3} that $\ox$ is fully stable in \eqref{5.30}, and thus the solution map $\breve S$ has a single-valued Lipschitz continuous localization $\vt$ around $(\ov,\op,\ox)$ with $\ov=(0,0,0)$. We also deduce from \cite[Example~6.4]{MN} that GSSOSC \eqref{6.8} fails when $\bar\lm:=(\frac{3}{8},\frac{5}{8},0,0)$.

Let us now check that the aforementioned SCOC does not hold here. Indeed, note that $\bar\lm$ is an extreme point of the set $\Lambda(\ox,\op,\ov)$, which can be calculated directly by \eqref{6.6} as
$$
\Lambda(\ox,\op,\ov)=\Big\{\Big(\frac{3}{8}-\al,\frac{5}{8}-\al,\al,\al\Big)\Big|\;0\le\al\le\frac{3}{8}\Big\}.
$$
Observe further that the gradient vectors $\nabla_x\ph_1(\ox,\op),\nabla_x\ph_2(\ox,\op)$ are linearly independent in $\R^3$ and that the determinant of the matrix
\[\begin{pmatrix}
\nabla_x L(\ox,\op,\bar\lm)&\nabla\ph_1(\ox)^T&\nabla\ph_2(\ox)^T\\
-\nabla\ph_1(\ox)&0&0\\
-\nabla\ph_2(\ox)&0&0
\end{pmatrix}=\begin{pmatrix}
0&-1&0&1&-1\\
-1&0&0&0&0\\
0&0&2&-1&-1\\
-1&0&1&0&0\\
1&0&1&0&0
\end{pmatrix}
\]
is equal to zero. It demonstrates the violation of SCOC from \cite[Definition~5.4.11]{FP} in this example.}
\end{Example}\vspace*{-0.07in}

Finally in this section, we discuss the role of the {\em pointbased} second-order condition GSSOSC in the study of full stability of PVC from \eqref{VC} and \eqref{6.2}. The following consequence of Theorem~\ref{thm6.3} describes the situation under the two first-order constraint qualifications considered above as well as under the partial {\em linear independence constraint qualification} (LICQ) at $(\ox,\op)$:
\[
\mbox{the gradients}\;\nabla_x\ph_i(\ox,\op)\;\mbox{ for }\;i\in I(\ox,\op)\;\mbox{ are linearly independent},
\]
which clearly implies both partial MFCQ and CRCQ at the corresponding point. Note that the assertions of the corollary below are new and rather surprising by taking into account the classical nature of the first-order and second-order conditions used and the novelty of the full stability notion for general PVC under consideration. On the other hand, we have recently obtained in \cite{MN} the prototypes of these results in parametric nonlinear programming.\vspace*{-0.1in}

\begin{Corollary} {\bf (full stability of PVC via GSSOSC).}\label{gssosc} Let $(\ox,\ov,\op,\oq)\in\R^n\times\R^n\times\R^l\times \R^m$ be such that $\ox\in\breve S(\ov,\op,\oq)$ in \eqref{6.3}, and let $f$ satisfy the standing assumptions of this section. Then the following assertions hold:

{\bf (i)} If both partial MFCQ and CRCQ are fulfilled at $(\ox,\op)$, then the validity of GSSOSC \eqref{6.8} ensures that $\ox$ is a fully stable solution of PVC in \eqref{VC} and \eqref{6.2} corresponding to $(\ov,\op)$.

{\bf (ii)} If the partial LICQ holds at $(\ox,\op)$, then GSSOSC is necessary and sufficient for the full stability of the corresponding solution $\ox$ in {\rm(i)}.
\end{Corollary}\vspace*{-0.07in}
{\bf Proof.} Since GSSOSC is stronger than GUSOSC as discussed after Definition~\ref{gusosc}, assertion (i) follows directly from Theorem~\ref{thm6.3}. To verify the necessary of GSSOSC for full stability in assertion (ii) under the partial LICQ, it suffices to show that the pointbased second-order subdifferential condition \eqref{4.43} from Proposition~\ref{Lipspq} is equivalent to GSSOSC. The proof of this fact follows the procedure in the proof of \cite[Theorem~6.6]{mrs}, where the inner product $\la z,w\ra$ with $(z,t)\in(D^*\partial_x g)(\ox,\op,\hat v)(w)$ and $\hat v=\ov-f(\ox,\op,\oq)$ is explicitly calculated.\endproof\vspace*{-0.15in}

\section{Concluding Remarks}
\setcounter{equation}{0}

This paper provides a systematical study of new H\"olderian and Lipschitzian stability notions for a general class of parametric variational systems as well as their significant specifications. There are many important issues remain for further research. Among them we list the following: deriving pointbased second-order conditions for H\"olderian full stability in general PVS and their PVI and PVC specifications in both finite and infinite dimensions under smoothness or/and strong monotonicity assumptions; developing the approach of this paper to establish efficient conditions for H\"olderian and Lipschitzian full stability of variational systems with nonpolyhedral conic constraints given, in particular, by in the parametric forms of second-order cone and semidefinite programming, etc. These and other open questions will be considered in our future work.
\end{document}